\documentclass[authoryear,1p]{elsarticle}

\usepackage{amsmath}
\usepackage{amssymb}
\usepackage{bm}
\usepackage{tikz}
\usepackage{pgfplots}
\usepackage{color, colortbl}
\usepackage{graphicx}

\graphicspath{{figs/}}

\usepackage{hyperref}
\usepackage{multirow}
\usepackage{float}

\newproof{pf}{Proof}

\journal{arXiv.org}

\begin{document}

\begin{frontmatter}

\title{Approximation of a fractional power of an elliptic operator}

\author[nsi,univ]{Petr N. Vabishchevich\corref{cor}}
\ead{vabishchevich@gmail.com}

\address[nsi]{Nuclear Safety Institute, Russian Academy of Sciences, 52, B. Tulskaya, Moscow, Russia}
\address[univ]{Peoples' Friendship University of Russia (RUDN University), 6 Miklukho-Maklaya St, Moscow, Russia}

\cortext[cor]{Corresponding author}

\begin{abstract}
Some mathematical models of applied problems lead to the need of solving
boundary value problems with a fractional power of an elliptic operator.
In a number of works, approximations of such a nonlocal operator are constructed on the basis of
an integral representation with a singular integrand.
In the present paper, new and more convenient integral representations are proposed
for operators with fractional powers. Approximations are based on the classical
quadrature formulas. The results of numerical experiments on the accuracy of
quadrature formulas are presented. The proposed approximations are used for numerical solving
a model two-dimensional boundary value problem for fractional diffusion.
\end{abstract}


\end{frontmatter}

\section{Introduction}\label{sec1}

Applied mathematical modeling of nonlocal processes \cite{baleanu2012fractional,uchaikin} is often associated with solving
boundary value problems that include fractional powers of elliptic operators \cite{Pozrikidis16}. 
For example, suppose that in a bounded domain $\Omega$ on the set of functions 
$v(\bm x) = 0, \ \bm x \in \partial \Omega$, 
there is defined the operator $\mathcal{A}$: $\mathcal{A}  v = - \triangle v, \ \bm x \in \Omega$. 
We seek the solution of the problem for the equation with the fractional power elliptic operator
$\mathcal{A}^{\alpha } v = f, \ 0 < \alpha < 1$ for a given $f(\bm x), \ \bm x \in \Omega$.

Various numerical techniques, such as finite difference or finite volume methods, can 
be used  to approximate problems with fractional power elliptic operators. 
The implementation of such algorithms require to calculate the action of a matrix (operator)
function on a vector $\varPhi(A) b$, where $A$ is a given 
matrix (operator) and $b$ is a given vector. For instance, in order to calculate 
the solution of the discrete fractional power elliptic problem, we get  
$\varPhi(z) = z^{-\alpha}$, where  $0 < \alpha < 1$.

There exist various approaches \cite{higham2008functions} how to calculate $\varPhi(A) b$.
They can be divided into two classes. In the first class, we construct one or another approximation of the operator
function, i.e., $R_m(A) \approx \varPhi^{-1}(A)$. In this case, the solution of the problem $R_m(A) u = b$ is connected, e.g.,
with the solution of a number of standard problems $(A + \gamma_i I) u_i = b_i, \ i = 1, \ldots, m$.
The second class of methods is based on the construction of approximate solutions $u \approx \varPhi(A) b$.

For a functional approximation of fractional powers, best uniform rational approximation \cite{stahl2003best} can be applied.
The use of such technology for approximate solving multidimensional problems of fractional diffusion is
discussed, for example, in \cite{harizanov2018optimal}. 

For a fractional power of a self-adjoint positive operator, the following integral
representation holds \cite{balakrishnan1960fractional,birman1987spectral}:
\begin{equation}\label{1}
 A^{-\alpha} = \frac{\sin(\alpha \pi)}{\pi} \int_{0}^{\infty } \theta^{-\alpha} (A + \theta I)^{-1} d \theta .
\end{equation}
To construct an approximation formula, the corresponding quadrature formulas \cite{bonito2015numerical} are used.
The fractional power of the operator is approximated by the sum of standard operators.
To improve the accuracy of approximation, various approaches are applied.

The kernel of the representation (\ref{1}) has a peculiarity if $\theta \rightarrow 0$.
To eliminate it, it is possible \cite{bonito2016numerical} to introduce the new variable $\theta = e^t$.
Thus, from (\ref{1}), we arrive at
\[
 A^{-\alpha} = \frac{\sin(\alpha \pi)}{\pi} \int_{-\infty }^{\infty } e^{(1-\alpha)t} (A + e^t I)^{-1} d t .
\]
To obtain the integral on a finite interval, in \cite{frommer2014efficient}, there is employed the transformation
\[
 \theta = \mu \frac{1-\eta}{1+\eta},
 \quad \mu > 0.
\] 
From (\ref{1}), we have
\begin{equation}\label{2}
  A^{-\alpha  } = \frac{2 \mu^{1-\alpha } \sin(\pi \alpha  )}{\pi} \int_{-1}^{1} (1-t)^{-\alpha }(1+t)^{\alpha -1}
  \big (\mu (1-t) I + (1+t) A \big )^{-1} d t . 
\end{equation} 
To approximate the right-hand side, we apply the Gauss-Jacobi quadrature formula \cite{Rabinowitz} 
with the weight $(1-t)^{-\alpha }(1+t)^{\alpha -1}$.

In the present paper, a more general integral representation (in comparison with (\ref{1}))
is used for fractional power operators, and possibilities of constructing new approximations are considered.
Section 2 describes a problem for an equation with a fractional power of an elliptic operator of second order.
Various types of integral representations of the fractional power operator are discussed in Section 3.
In Section 4, we construct the corresponding quadrature formulas.
Examples of solving a model two-dimensional problem with the fractional
power of the Laplace operator are given in Section 5.
The results of the work are summarized in Section 6.

\section{Problem formulation}\label{sec:2}

In a bounded polygonal domain $\Omega \subset R^2$ with the Lipschitz continuous boundary $\partial\Omega$,
we search the solution of a  problem with a fractional power of an elliptic operator.
Here, we use the definition of a fractional power of an elliptic operator that relies on the spectral 
theory \cite{birman1987spectral,kwasnicki2017ten}.
Introduce the elliptic operator as
\begin{equation}\label{3}
  \mathcal{A}  v = - {\rm div}  ( a({\bm x}) {\rm grad} \, v) + c({\bm x}) v
\end{equation} 
with coefficients $0 < a_1 \leq a({\bm x}) \leq a_2$, $c({\bm x}) \geq 0$.
The operator $\mathcal{A}$ is defined on the set of functions $v({\bm x})$ that satisfy
on the boundary $\partial\Omega$ the following conditions:
\begin{equation}\label{4}
  v ({\bm x}) = 0,
  \quad {\bm x} \in \partial \Omega .
\end{equation} 

In the Hilbert space $\mathcal{H} = L_2(\Omega)$, we define 
the scalar product and norm in the standard way:
\[
  (v,w) = \int_{\Omega} v({\bm x}) w({\bm x}) d{\bm x},
  \quad \|v\| = (v,v)^{1/2} .
\] 
In the spectral problem
\[
 \mathcal{A}  \varphi_k = \lambda_k \varphi_k, 
 \quad \bm x \in \Omega , 
\] 
\[
  \varphi_k = 0,
  \quad {\bm x} \in \partial \Omega , 
\] 
we have 
\[
 0 < \lambda_1 \leq \lambda_2 \leq ... ,
\] 
and the eigenfunctions  $ \varphi_k, \ \|\varphi_k\| = 1, \ k = 1,2, ...  $ form a basis in $L_2(\Omega)$. Therefore, 
\[
 v = \sum_{k=1}^{\infty} (v,\varphi_k) \varphi_k .
\] 
Let the operator $\mathcal{A}$ be defined in the following domain:
\[
 D(\mathcal{A} ) = \Big \{ v \ | \ v(\bm x) \in L_2(\Omega), \ \sum_{k=0}^{\infty} | (v,\varphi_k) |^2 \lambda_k < \infty \Big \} .
\] 
Under these conditions  $\mathcal{A} : L_2(\Omega) \rightarrow L_2(\Omega)$ and
the operator $\mathcal{A}$ is self-adjoint and positive definite: 
\begin{equation}\label{5}
  \mathcal{A}  = \mathcal{A} ^* \geq \nu   \mathcal{I} ,
  \quad \nu   > 0 ,    
\end{equation} 
where $\mathcal{I}$ is the identity operator in $\mathcal{H}$.
In applications, the value of $\lambda_1$ is unknown (the spectral problem must be solved).
Therefore, we suppose that $\nu  \leq \lambda_1$ in (\ref{5}).
Let us assume for the fractional power of the  operator $\mathcal{A}$:
\[
 \mathcal{A} ^\alpha v =  \sum_{k=0}^{\infty} (v,\varphi_k) \lambda_k^\alpha  \varphi_k .
\] 
More general and mathematically complete definition of fractional powers of elliptic operators 
is given in the work \cite{carracedo2001theory,yagi2009abstract}. 
The solution $v(\bm x)$ satisfies the equation
\begin{equation}\label{6}
  \mathcal{A}^\alpha v = f 
\end{equation} 
under the restriction $0 < \alpha < 1$. 

We consider the simplest case, where the computational domain $\Omega$ is a rectangle
\[
 \Omega = \{ \bm x  \ | \ \bm x = (x_1,x_2), \ 0 < x_n < l_n, \ n = 1,2 \} .
\]
To solve approximately the problem (\ref{6}), we introduce in the domain $\Omega$ a uniform grid
\[
\overline{\omega}  = \{ \bm{x} \ | \ \bm{x} =\left(x_1, x_2\right), \quad x_n =
i_n h_n, \quad i_n = 0,1,...,N_n,
\quad N_n h_n = l_n, \ n = 1,2 \} ,
\]
where $\overline{\omega} = \omega \cup \partial \omega$ and
$\omega$ is the set of interior nodes, whereas $\partial \omega$ is the set of boundary nodes of the grid.
For grid functions $u(\bm x)$ such that $u(\bm x) = 0, \ \bm x \notin \omega$, we define the Hilbert space
$H=L_2\left(\omega\right)$, where the scalar product and the norm are specified as follows:
\[
\left(u, w\right) =  \sum_{\bm x \in  \omega} u\left(\bm{x}\right)
w\left(\bm{x}\right) h_1 h_2,  \quad 
\| y \| =  \left(y, y\right)^{1/2}.
\]

For $u(\bm x) = 0, \ \bm x \notin \omega$, the grid operator $A$  can be written as
\[
  \begin{split}
  A u = & -
  \frac{1}{h_1^2} a(x_1+0.5h_1,x_2) (u(x_1+h_1,h_2) - u(\bm{x})) \\ 
  & + \frac{1}{h_1^2} a(x_1-0.5h_1,x_2) (u(\bm{x}) - u(x_1-h_1,h_2)) \\
  & - \frac{1}{h_2^2} a(x_1,x_2+0.5h_2) (u(x_1,x_2+h_2) - u(\bm{x})) \\ 
  & + \frac{1}{h_2^2} a(x_1, x_2-0.5h_2) (u(\bm{x}) - u(x_1,x_2-h_2)) + c(\bm x) u(\bm x), 
  \quad \bm{x} \in \omega . 
 \end{split} 
\] 
For the above grid operators (see \cite{Samarskii1989,SamarskiiNikolaev1978}), we have 
\[ A = A^* \geq \delta I,
 \quad \delta > 0 ,
\]
where $I$ is the grid identity operator. 
For problems with sufficiently smooth coefficients and the right-hand side, 
it approximates the differential operator with the truncation error 
$\mathcal{O} \left(|h|^2\right)$, $|h|^2 = h_1^2+h_2^2$. 

To handle the fractional power of the grid operator $A$, let us consider the eigenvalue problem
\[
 A \psi_k = \mu_k \psi_k . 
\] 
We have
\[
 0 < \delta = \mu_1 \leq \mu_2 \leq ... \leq \mu_K,
 \quad K = (N_1-1)(N_2-1) , 
\] 
where eigenfunctions $\psi_k, \ \|\psi_k\| = 1, \ k = 1,2, ..., K,$ form a basis in $H$. Therefore
\[
 u = \sum_{k= 1}^{K}(u, \psi_k) \psi_k . 
\]
For the fractional power of the operator $A$, we have
\[
 A^\alpha u = \sum_{k= 1}^{K}(u, \psi_k) \mu_k^\alpha \psi_k .
\] 

Using the above approximations, from (\ref{6}), we arrive at the discrete problem
\begin{equation}\label{7}
 u = A^{-\alpha} b .
\end{equation} 
Taking into account the transformation
\[
 A \longrightarrow \frac{1}{\delta } A,
 \quad b \longrightarrow \frac{1}{\delta^\alpha } b ,
\]
we can assume, without loss of generality, that in (\ref{7}), we have
\begin{equation}\label{8}
 A \geq I,
\end{equation} 
i.e., $\delta = 1$.

\section{Integral representation of the fractional power operator}\label{sec:3}

For numerical solving the problem (\ref{7}), (\ref{8}), it is necessary to have an integral representation of the function
$x^{-\alpha}$ for $x > 1$ and $0 < \alpha < 1$.
The starting point of our constructions is the formula (see the definite integral 3.141.4 in the book \cite{gradshteyn2007table})
\begin{equation}\label{9}
 \int_{0}^{\infty} \frac{\theta^{\mu-1}}{(p + q \theta^\nu)^{n+1}} d \theta = q^{-\mu /\nu} p^{\mu /\nu - n-1} \gamma \Big (n, \nu, \frac{\mu }{\nu } \Big), 
\end{equation} 
where
\[
 0 < \frac{\mu }{\nu } < n+1,
 \quad p \neq 0,  \quad q \neq 0,
 \quad \gamma(n, \nu, \xi ) = \frac{1}{\nu } \frac{\Gamma (\xi) \, \Gamma (1+n-\xi)}{\Gamma (1+n)} ,
\] 
and $\Gamma(x)$ is the gamma function.
Using (\ref{9}), we can construct different integral representations for the fractional power operator.
Let us discuss some possibilities in this direction.

To employ  (\ref{9}), we introduce parameters $p, q$ associated with the operator $A$:  $p = p(A), \ q=q(A)$.
In choosing these dependencies, we must bear in mind that it is necessary to calculate $(p(A) + q(A) \theta^\nu)^{-1}$  for
integers $n \geq 0$. With this in mind, it is natural to put $n = 0$ and consider two possibilities:
\begin{equation}\label{10}
 p \rightarrow  A,
 \quad q \rightarrow  I, 
\end{equation} 
\begin{equation}\label{11}
 p\rightarrow  I,
 \quad q \rightarrow A . 
\end{equation} 

Using (\ref{10}), from (\ref{9}), we get the following integral representation of the fractional power operator:
\begin{equation}\label{12}
 A^{\mu /\nu - 1} = \frac{\nu }{\pi } \sin \Big (\pi \frac{\mu }{\nu } \Big ) 
 \int_{0}^{\infty} \theta^{\mu-1} (A + \theta^\nu I)^{-1} d \theta .
\end{equation}
In the particular case of $\nu = 1, \ \mu = 1 - \alpha$, from (\ref{12}), it follows the Balakrishnan formula (\ref{1}).
The inverse transition also takes place, i.e., when using the new variable $\theta = t^\nu$, from (\ref{1}), we arrive at (\ref{12}).
For (\ref{12}), from (\ref{9}), we obtain the representation
\begin{equation}\label{13}
 A^{-\mu /\nu} = \frac{\nu }{\pi } \sin \Big (\pi \frac{\mu }{\nu } \Big ) 
 \int_{0}^{\infty} \theta^{\mu-1} ( I + \theta^\nu A)^{-1} d \theta .
\end{equation}
This representation corresponds to the use of the new variable $\theta = t^{-\nu}$ in (\ref{1}).
Thus, for $n = 0$, the integral representations (\ref{12}), (\ref{13}) can be treated as variants of the
Balakrishnan formula (\ref{1}). Obviously, for $n > 0$, we do not have such the direct relation between (\ref{1})
and (\ref{9}), (\ref{10}) or (\ref{9}), (\ref{11}).

In the two-parameter formulas (\ref{12}) and (\ref{13}), one parameter is free.
If we put $\mu = 1$, then the integrand in the right-hand side of (\ref{12}) and (\ref{13})
have no singularities at the ends of the integration interval. Assuming $\alpha = 1 - 1 /\nu$,
from (\ref{12}), we get
\begin{equation}\label{14}
 A^{-\alpha} = \frac{\sin (\pi \alpha ) }{(1-\alpha) \pi }
 \int_{0}^{\infty} (A + \theta^{1/(1-\alpha)} I)^{-1} d \theta .
\end{equation} 
In case of (\ref{13}), we have $\alpha = 1/\nu$ and
\begin{equation}\label{15}
 A^{-\alpha} = \frac{\sin (\pi \alpha )}{\alpha \pi }  
 \int_{0}^{\infty} (I + \theta^{1/\alpha} A)^{-1} d \theta .
\end{equation} 

The ability to use standard quadrature formulas for approximating the fractional power operator
is provided by passing from (\ref{14}) and (\ref{15}) to a finite integral.
Let us consider, e.g., the transition to a new integration variable $t$ on the interval $[0,1]$ such
that $\theta = t (1-t)^\sigma$ with a numerical parameter $\sigma > 0$.
For (\ref{14}), we get
\begin{equation}\label{16}
 A^{-\alpha} = \frac{\sin (\pi \alpha )}{(1-\alpha) \pi } 
 \int_{0}^{1} (1-t)^{\sigma \alpha /(1-\alpha) - 1}  \big (1+(\sigma -1)t \big ) \big((1-t)^{\sigma /(1-\alpha)} A + t^{1/(1-\alpha)} I \big)^{-1} d t .
\end{equation}

The integrand in (\ref{16}) has no singularities for $\sigma \geq \sigma_1 = (1-\alpha)/\alpha$.
For  $\sigma =  \sigma_1$, we have
\begin{equation}\label{17}
 A^{-\alpha} = \frac{\sin (\pi \alpha )}{\alpha (1-\alpha) \pi } 
 \int_{0}^{1} \big(\alpha +(1-2\alpha)t \big) \big((1-t)^{1 /\alpha} A + t^{1/(1-\alpha)} I \big)^{-1} d t .
\end{equation}

The same change of the integration variable in the representation (\ref{15}) results in
\begin{equation}\label{18}
 A^{-\alpha} = \frac{\sin (\pi \alpha )}{\alpha \pi } 
 \int_{0}^{1} (1-t)^{\sigma (1-\alpha) /\alpha - 1}  
\big(1+(\sigma -1)t \big) \big((1-t)^{\sigma /\alpha} I + t^{1/\alpha} A \big)^{-1} d t . 
\end{equation}
We use this representation for $\sigma \geq \sigma_2$, where $\sigma_2 = \alpha / (1-\alpha)$.
For $\sigma = \sigma_2$, we get
\begin{equation}\label{19}
 A^{-\alpha} = \frac{\sin (\pi \alpha )}{\alpha (1-\alpha) \pi } 
 \int_{0}^{1} \big(1-\alpha  + (2\alpha -1)t \big) \big((1-t)^{1 /(1-\alpha)} I + t^{1/\alpha} A \big)^{-1} d t .
\end{equation}
Within the replacement of $t$ by $1-t$, the formulas (\ref{17}) and (\ref{19}) coincide.
The appropriate choice of $\sigma$ allows us to provide the $l$-th derivative of the integrand
without any singularities. In particular, for (\ref{16}) and (\ref{18}), it is necessary to take
$\sigma \geq  \sigma_1(l+1)$ and $\sigma \geq  \sigma_2(l+1)$, respectively.
This control of derivatives of the integrand is necessary when constructing quadrature formulas.

It seems reasonable to construct integral representations of the fractional power operator, when we select immediately
the integration interval from 0 to 1. We can use, e.g., the definite integral 3.197.10 from the book \cite{gradshteyn2007table}:
\[
 \int_{0}^{1} \frac{\theta^{\mu -1} (1-\theta)^{-\mu}} {1 + p \theta} d \theta = (1+p)^{-\mu} \frac{\pi }{\sin(\pi \mu)} ,
 \quad 0 < \mu < 1, \quad p > - 1.  
\]
Using the association $p \rightarrow A-I$, we arrive at the following integral representation:
\begin{equation}\label{20}
  A^{-\alpha} = \frac{\sin (\pi \alpha )}{\pi } 
 \int_{0}^{1} \theta^{\alpha -1} (1-\theta)^{-\alpha } \big( I + \theta (A-I)\big)^{-1} d t .
\end{equation}
Introducing the new variable $\theta = (1+t)/2$, we see that (\ref{20}) corresponds to the representation (\ref{2}) with $\mu = 1$.

The following formula (see, for example, 3.198 in the book \cite{gradshteyn2007table}) provides more possibilities:
\[
 \int_{0}^{1} \frac{\theta^{\mu -1} (1-\theta)^{\nu -1}} {( a \theta  + b(1-\theta + c)^{\mu +\nu}} d \theta = (a+c)^{-\mu} (b+c)^{-\nu } 
 \mathrm{B}   (\mu , \nu),
\]
\[
 \quad a \geq 0, \quad b \geq 0, \quad c > 0,
 \quad \mu > 0, \quad \nu > 0,
\] 
where $\mathrm{B}  (\mu , \nu)$ is the beta function. This representation can be used for
positive integers $n = \mu +\nu$. However, for $n=1$ and even $n=2$, the integrand has a singularity.

We also highlight the formula 3.234.2 in the book \cite{gradshteyn2007table}:
\[
 \int_{0}^{1} \left (
 \frac{\theta^{\mu -1}}{1 + p \theta} + \frac{\theta^{-\mu}}{p+\theta}\right )  d \theta = p^{-\mu} \frac{\pi }{\sin(\pi \mu)} ,
 \quad 0 < \mu < 1, \quad p > 0.  
\]
It is also inconvenient for the direct application due to the singularity of the integrand for $\theta=0$.
Using the new variable $\theta = t^\sigma$, we arrive at
\[
 \int_{0}^{1} \left (
 \frac{t^{\sigma \mu -1}}{1 + p t^\sigma } + \frac{t^{\sigma(1-\mu)-1}}{p+t^\sigma }\right )  d t = p^{-\mu} \frac{\pi }{\sigma \sin(\pi \mu)} . 
\]
The integrand has no sigularities if
\[
 \sigma \mu -1 \geq  0,
 \quad  \sigma(1-\mu)-1 \geq 0 .
\]
For instance, it is sufficient to take $\sigma \geq \sigma_3$, where $\sigma_3 = \max\big(\mu^{-1}, (1-\mu)^{-1}\big)$.
With this in mind, we obtain a new integral representation of the fractional power operator:
\begin{equation}\label{21}
 A^{-\alpha} = \frac{\sigma \sin (\pi \alpha )}{\pi }  \int_{0}^{1}  \Big(t^{\sigma \alpha  -1} \big(I +  t^\sigma A\big)^{-1} 
 + t^{\sigma(1-\alpha )-1} \big(A+t^\sigma I \big)^{-1} \Big ) d t .
\end{equation}
Selecting $\sigma \geq  \sigma_3(l+1)$, we guarantee that the $l$-th derivative of the integrand has no singularities.

The integral representations (\ref{17}) (or (\ref{19})) and (\ref{21}) are applied to approximate the
fractional power operator using one or another quadrature formulas.

\section{Numerical integration}\label{sec:4} 

The accuracy of quadrature formulas is investigated for the following integral
representation (see (\ref{16})) of the function $x^{-\alpha}$ with $x \geq 1$:
\begin{equation}\label{22}
\begin{split}
 x^{-\alpha} & = \frac{\sin (\pi \alpha )}{(1-\alpha) \pi } 
 \int_{0}^{1} (1-t)^{\sigma \alpha /(1-\alpha) - 1}  \big (1+(\sigma -1)t \big ) 
\big((1-t)^{\sigma /(1-\alpha)} x + t^{1/(1-\alpha)}  \big)^{-1} d t , \\
 \quad \sigma & = \varkappa \frac{1-\alpha }{\alpha } , 
\end{split}
\end{equation} 
with the parameter $\varkappa \geq 1$.
When focusing on (\ref{21}), we are interested in the integral representation
\begin{equation}\label{23}
\begin{split}
 x^{-\alpha} & = \frac{\sigma \sin (\pi \alpha )}{\pi }  \int_{0}^{1}  \Big(t^{\sigma \alpha  -1} \big(1+  t^\sigma x\big)^{-1} 
 + t^{\sigma(1-\alpha )-1} \big(x+t^\sigma \big)^{-1} \Big ) d t , \\
 \quad \sigma & = \varkappa \max\big(\alpha^{-1}, (1-\alpha)^{-1}\big) .
\end{split}
\end{equation} 

\begin{figure}
\centering
\begin{minipage}{0.43\linewidth}
\centering
\includegraphics[width=\linewidth]{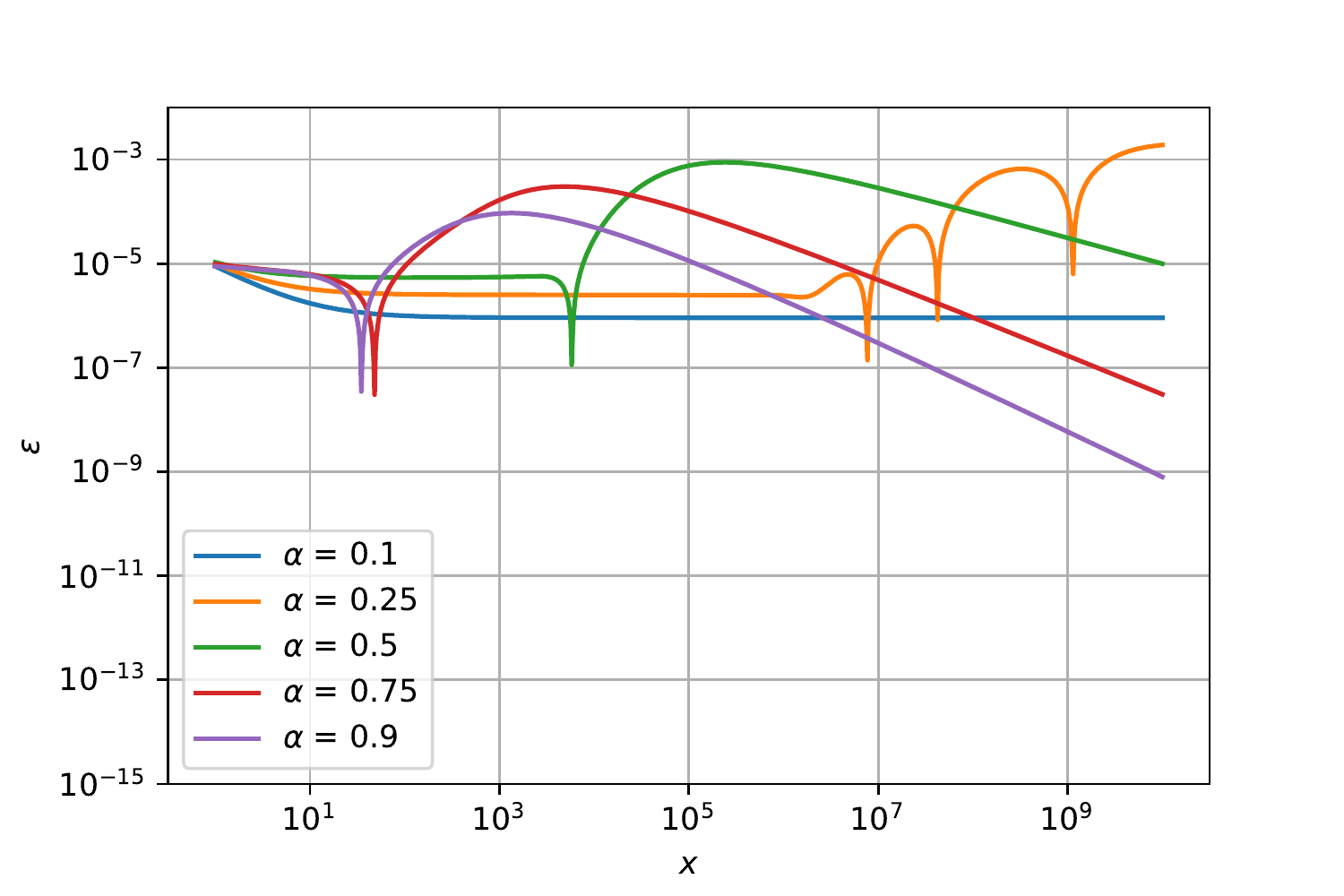}\\
$\varkappa = 1$ \\
\includegraphics[width=\linewidth]{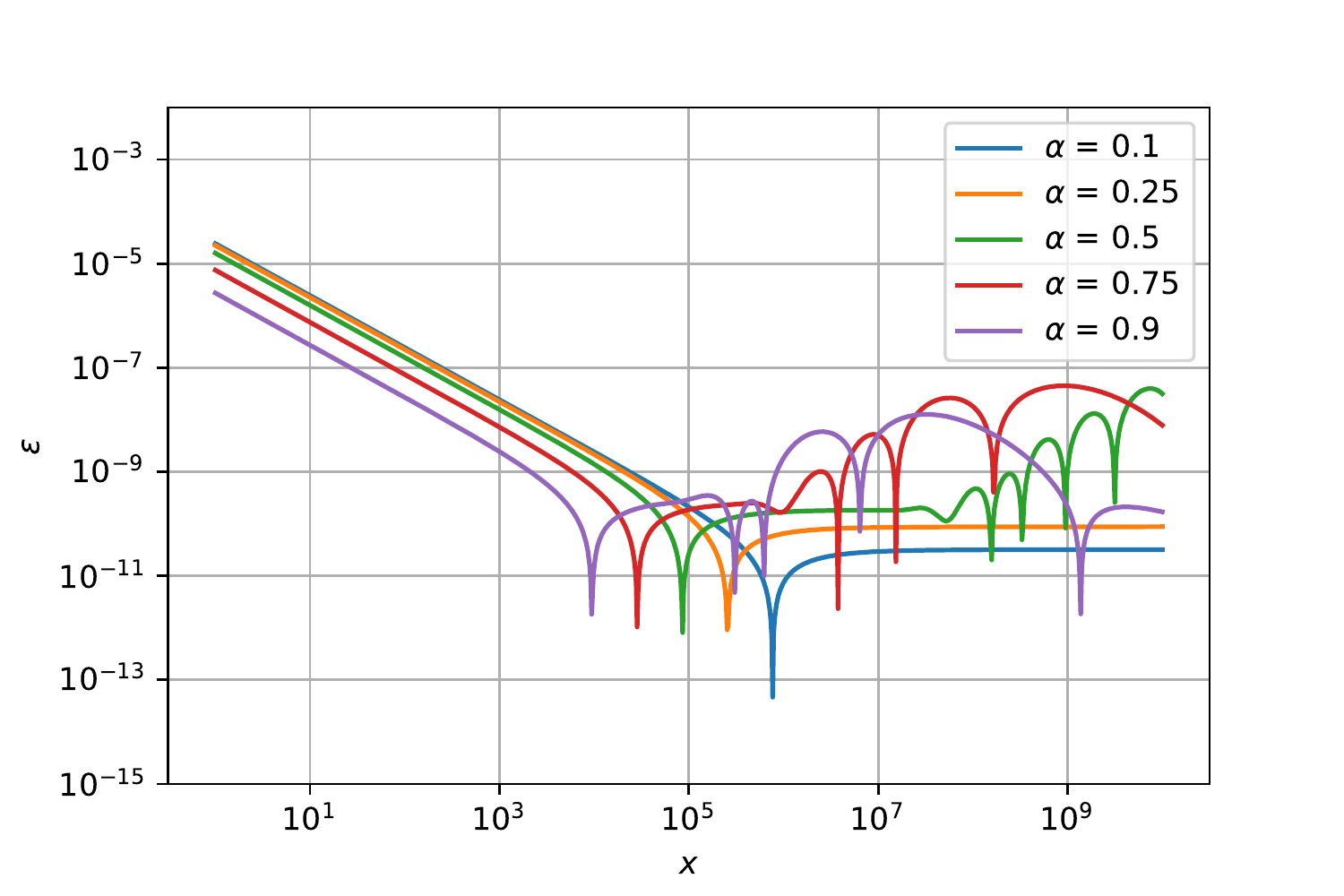}\\
$\varkappa = 3$ \\
\includegraphics[width=\linewidth]{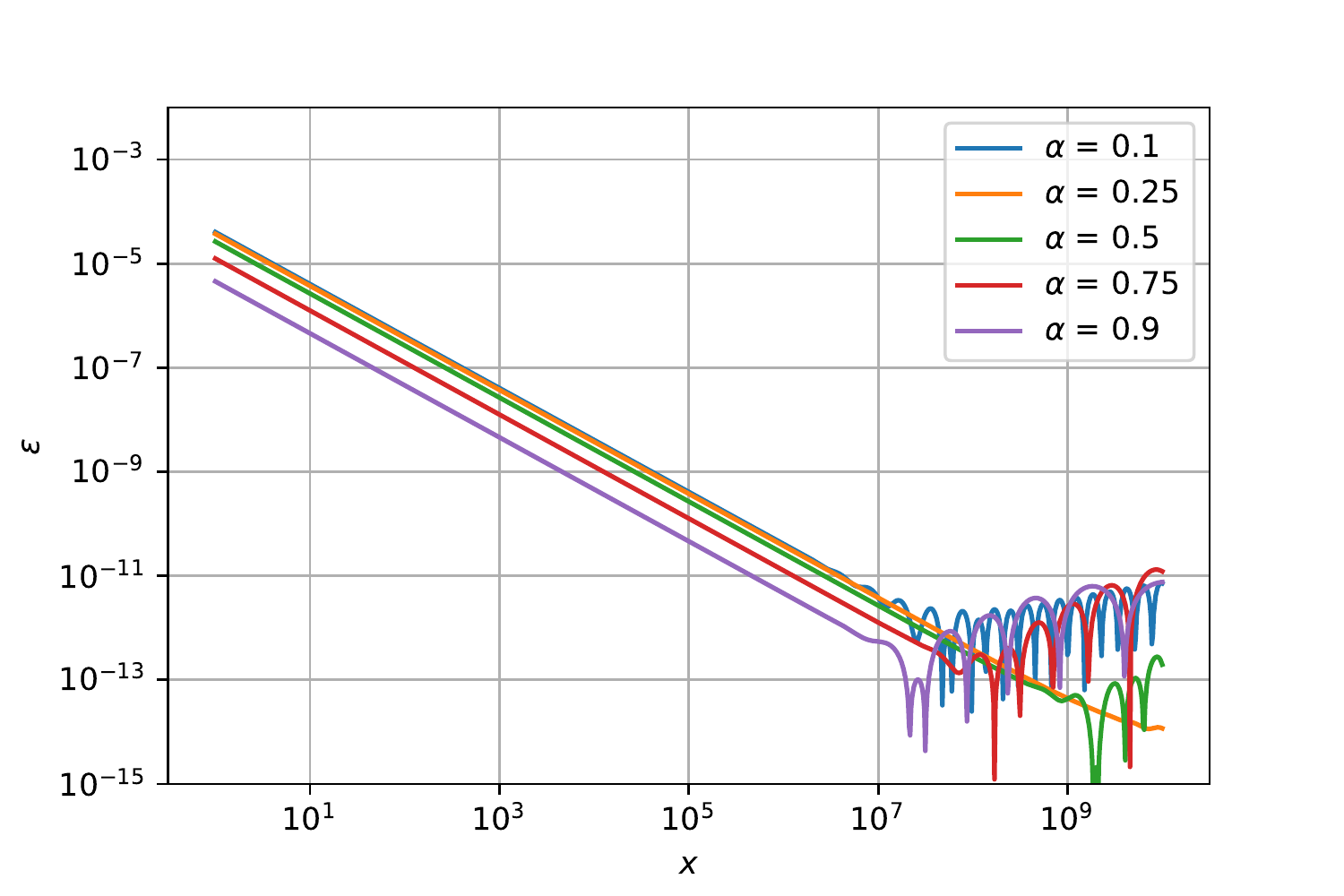}\\
$\varkappa = 5$ \\
\end{minipage}
\begin{minipage}{0.43\linewidth}
\centering
\includegraphics[width=\linewidth]{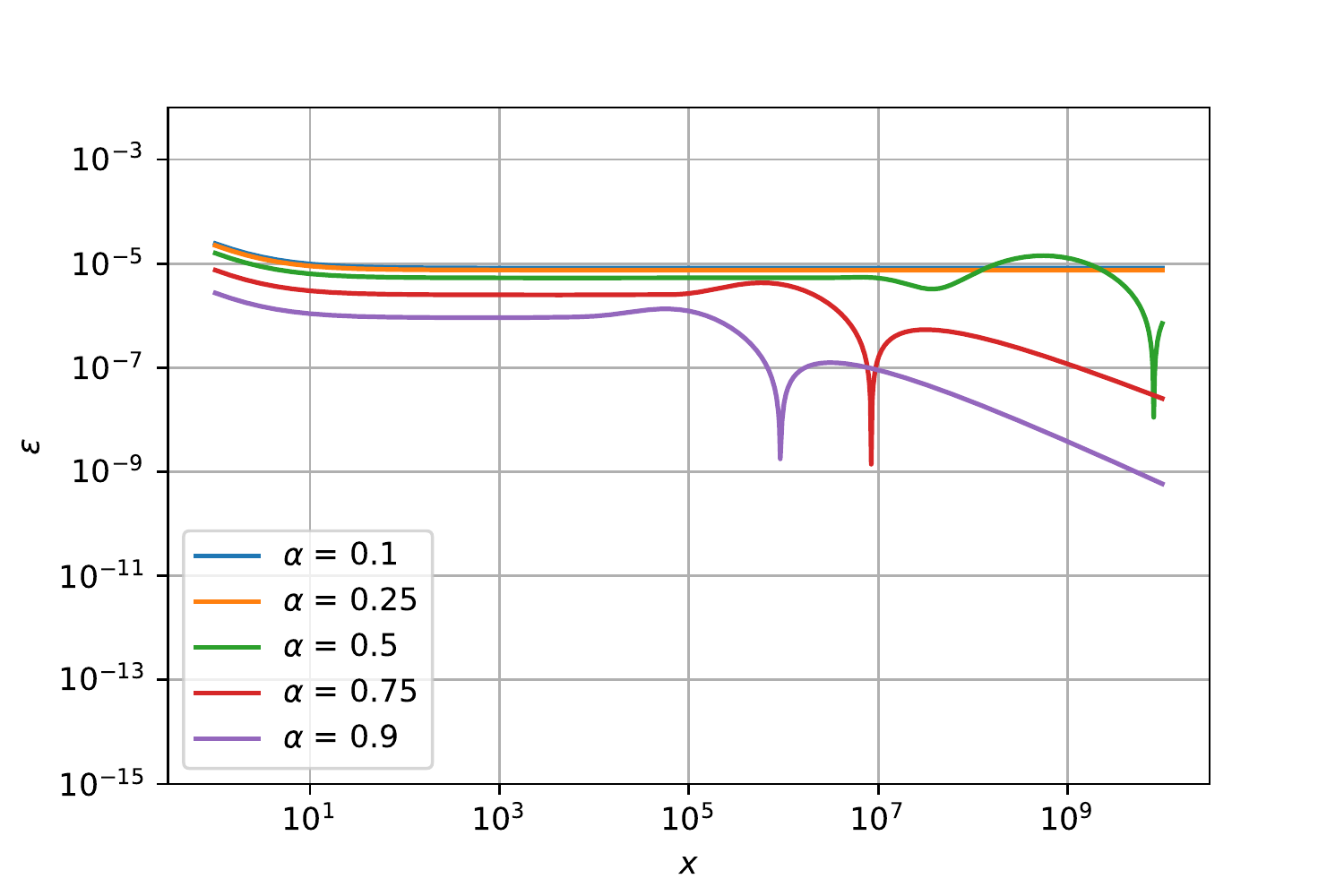}\\
$\varkappa = 2$ \\
\includegraphics[width=\linewidth]{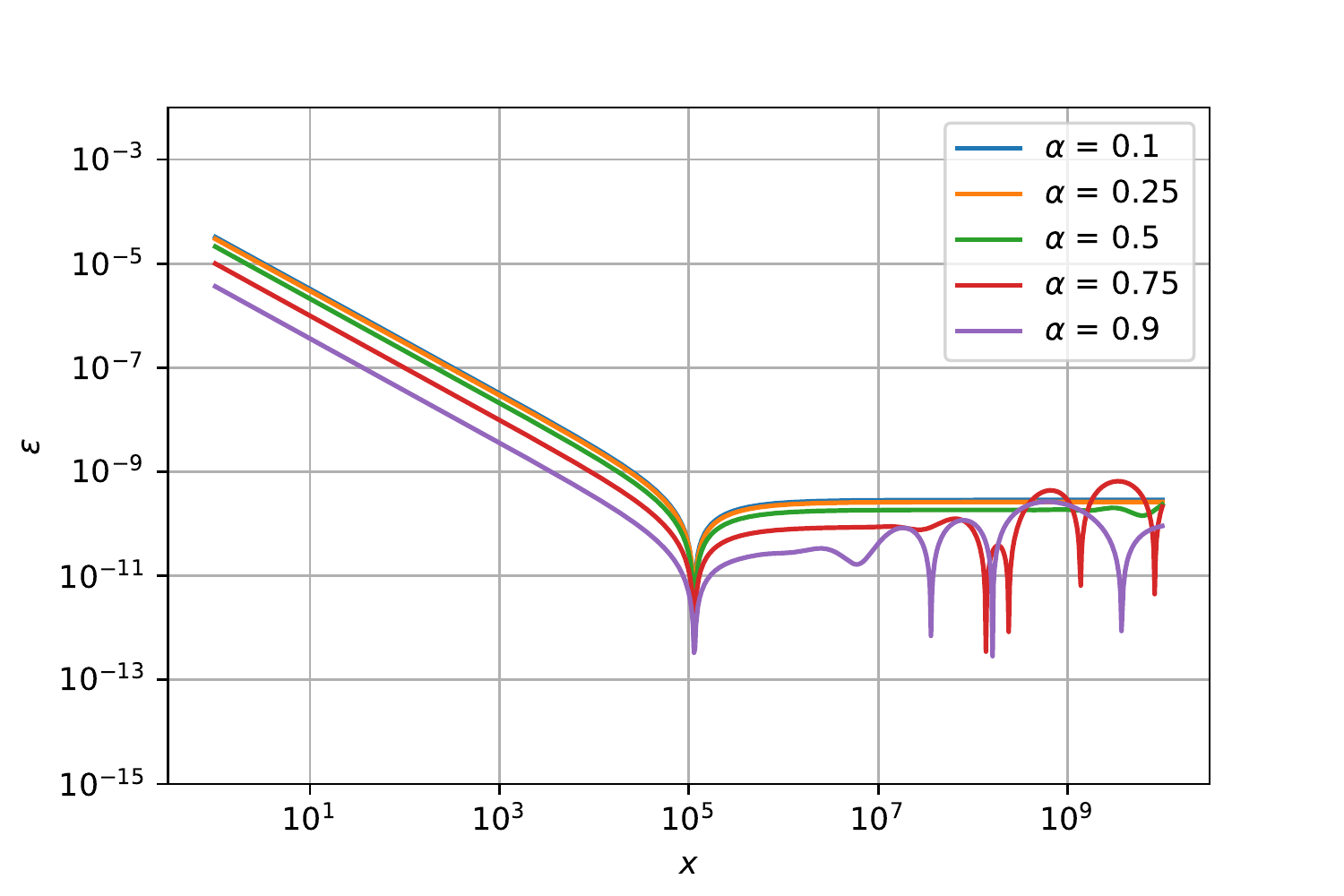}\\
$\varkappa = 4$ \\
\includegraphics[width=\linewidth]{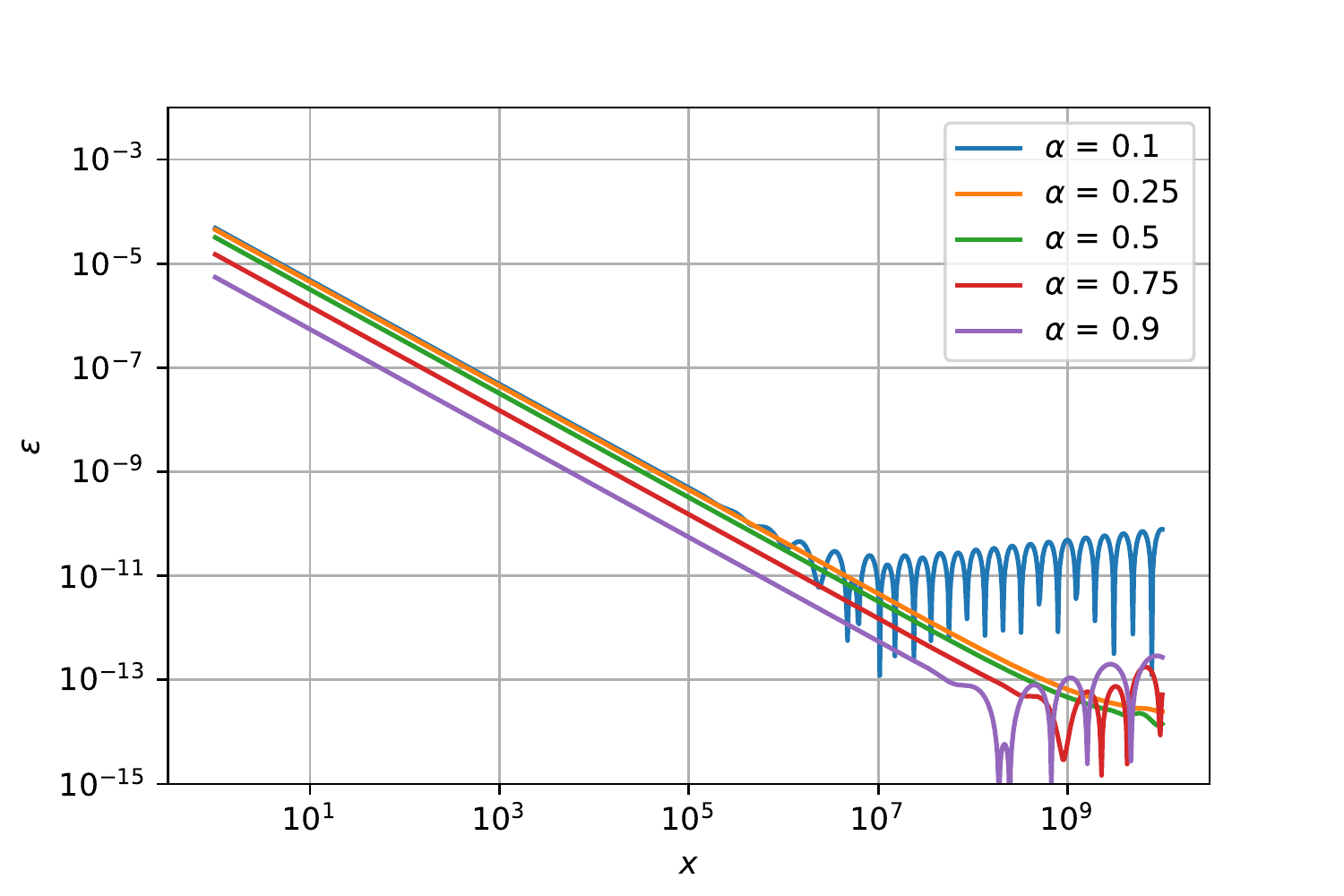}\\
$\varkappa = 6$ \\
\end{minipage}
\caption{Approximation error: the rectangle rule for (\ref{22}), $M = 100$.}
\label{f-1}
\end{figure}

We apply the standard quadrature formulas, i.e., the rectangle and Simpson's quadrature rules \cite{DavisRabinowitz}. 
The integration interval $[0,1]$ is divided into $M$ equal parts. When using the rectangle rule,
the integrand is calculated at the midpoint of these parts, whereas Simpson's rule employs
their boundary nodes. For (\ref{22}) and (\ref{23}) with the rectangle rule, the approximating function is denoted by
$r_1(x,\alpha; M, \varkappa)$ and $r_2(x,\alpha; M, \varkappa)$, respectively. 
Here $M$ and $\varkappa$ are explicitly highlighted as the key approximation parameters.
The approximation error is estimated at $1 \leq x \leq 10^{10}$ and is evaluated as follows:
\[
 \varepsilon_e(x) = |r_e(x,\alpha; M, \varkappa) - x^{-\alpha}|,
 \quad e = 1,2 .   
\]
For the rectangle rule, the approximation error of the function $ x^{-\alpha}$ 
is shown in Figure~\ref{f-1} for $\alpha = 0.1, 0.25, 0.5, 075, 0.9$
and various values of the parameter $\varkappa$. The calculations were performed at $M = 100$.
The significant effect of the parameter $\varkappa$ is observed, when $\varkappa = 1$ and $2$. 
The accuracy of the rectangle quadrature formula is limited by the second derivative of the integrand.
Therefore, we must take $\varkappa \geq 2$. To improve the approximation accuracy of the function $x^{-\alpha}$ for large
values of $x$, we need to take larger values of the parameter $\varkappa$ .

Numerical data for the maximal error
\[
 \overline{\varepsilon}_e = \max_{1 \leq x \leq 10^{20}} \varepsilon_e(x) ,
 \quad e = 1,2 ,  
\]
are presented in Table~1 for various numbers of integration nodes. It is easy to see that
when the number of partitions $M$ increases, the convergence rate becomes close to theoretical estimates.

\begin{center}
\begin{table}[htp]
\label{t-1}
\caption{Approximation error $\overline{\varepsilon}_e$: the rectangle rule for (\ref{22})}
\centering
\begin{tabular}{ccccccc}
\hline
   $M$   &    $\varkappa$   &   $\alpha = 0.1$         &     $\alpha = 0.25$   &    $\alpha = 0.5$   &     $\alpha = 0.75$      &    $\alpha = 0.9$   \\
\hline
        &      1   &  5.374386e-03     &  3.975611e-03     &    1.777798e-03    &  6.057090e-04    &  1.860264e-04  \\
        &      2   &  9.694686e-05     &  1.546665e-04     &    6.365306e-05    &  3.000009e-05    &  1.090598e-05  \\
   50   &      3   &  9.583340e-05     &  8.914710e-05     &    6.363818e-05    &  3.000215e-05    &  1.092388e-05  \\
        &      4   &  1.264541e-04     &  1.189327e-04     &    8.483208e-05    &  4.000540e-05    &  1.457062e-05  \\
        &      5   &  1.548898e-04     &  1.485671e-04     &    1.060141e-04    &  5.000617e-05    &  1.821371e-05  \\
        &      6   &  1.841219e-04     &  1.780018e-04     &    1.271753e-04    &  6.000429e-05    &  2.185621e-05  \\

\hline
        &      1   &  1.294685e-03     &  1.985566e-03     &    8.876394e-04    &  3.029813e-04    &  9.380258e-05  \\
        &      2   &  2.433513e-05     &  3.853785e-05     &    1.591494e-05    &  7.501982e-06    &  2.730852e-06  \\
  100   &      3   &  2.426684e-05     &  2.234929e-05     &    1.591401e-05    &  7.502111e-06    &  2.731974e-06  \\
        &      4   &  3.232626e-05     &  2.983899e-05     &    2.121750e-05    &  1.000297e-05    &  3.642973e-06  \\
        &      5   &  4.029470e-05     &  3.731919e-05     &    2.652025e-05    &  1.250368e-05    &  4.553745e-06  \\
        &      6   &  4.814931e-05     &  4.478680e-05     &    3.182170e-05    &  1.500422e-05    &  5.464479e-06  \\
\hline
        &      1   &  1.647562e-04     &  9.922321e-04     &    4.435056e-04    &  1.515251e-04    &  4.710344e-05  \\
        &      2   &  6.092649e-06     &  9.617959e-06     &    3.978839e-06    &  1.875618e-06    &  6.829863e-07  \\
  200   &      3   &  6.088573e-06     &  5.596207e-06     &    3.978781e-06    &  1.875626e-06    &  6.830565e-07  \\
        &      4   &  8.129377e-06     &  7.471178e-06     &    5.304967e-06    &  2.500845e-06    &  9.107634e-07  \\
        &      5   &  1.016443e-05     &  9.345551e-06     &    6.631108e-06    &  3.126054e-06    &  1.138456e-06  \\
        &      6   &  1.219226e-05     &  1.121914e-05     &    7.957167e-06    &  3.751252e-06    &  1.366146e-06  \\
\hline
\end{tabular}
\end{table}
\end{center}

\begin{figure}
\centering
\begin{minipage}{0.43\linewidth}
\centering
\includegraphics[width=\linewidth]{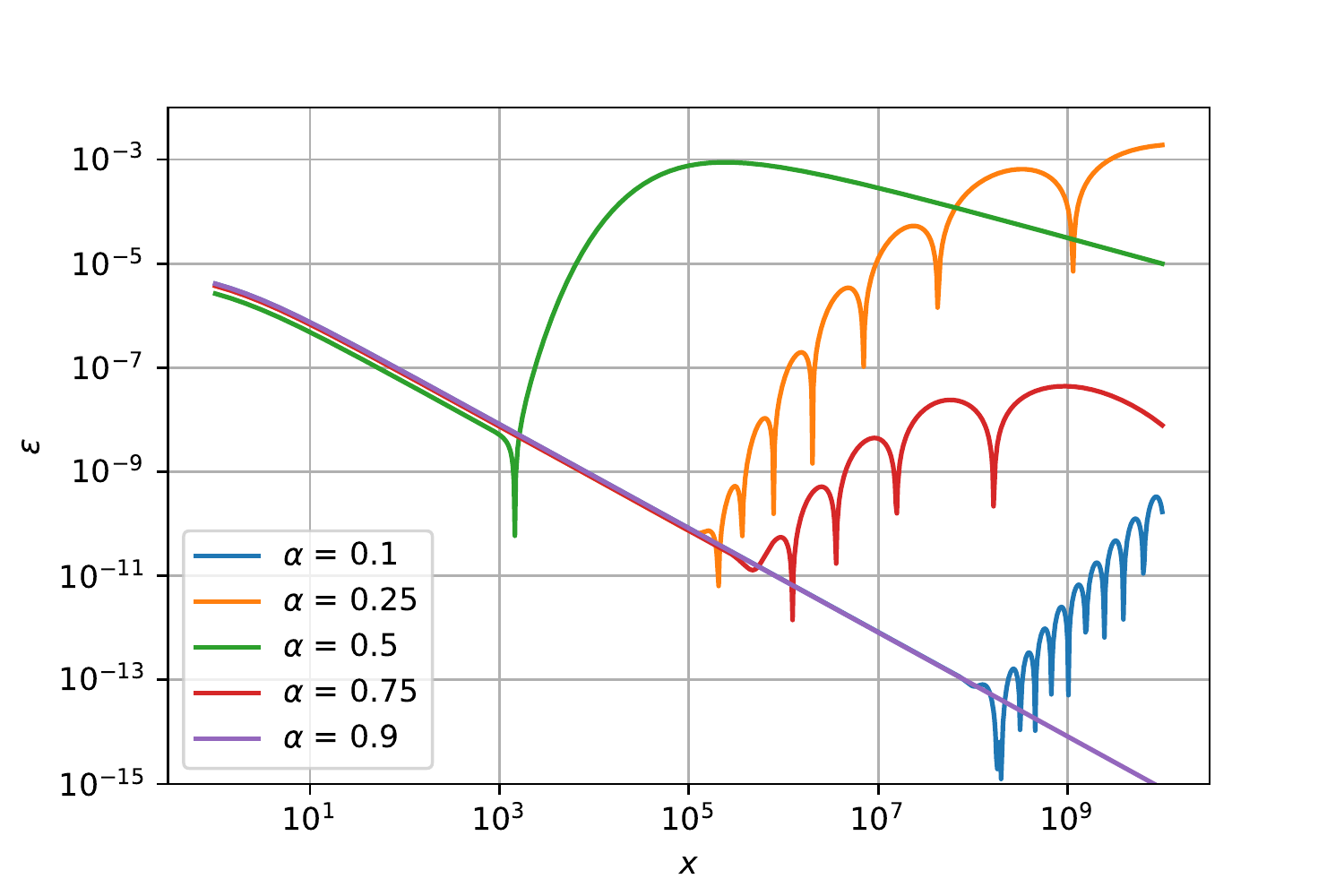}\\
$\varkappa = 1$ \\
\includegraphics[width=\linewidth]{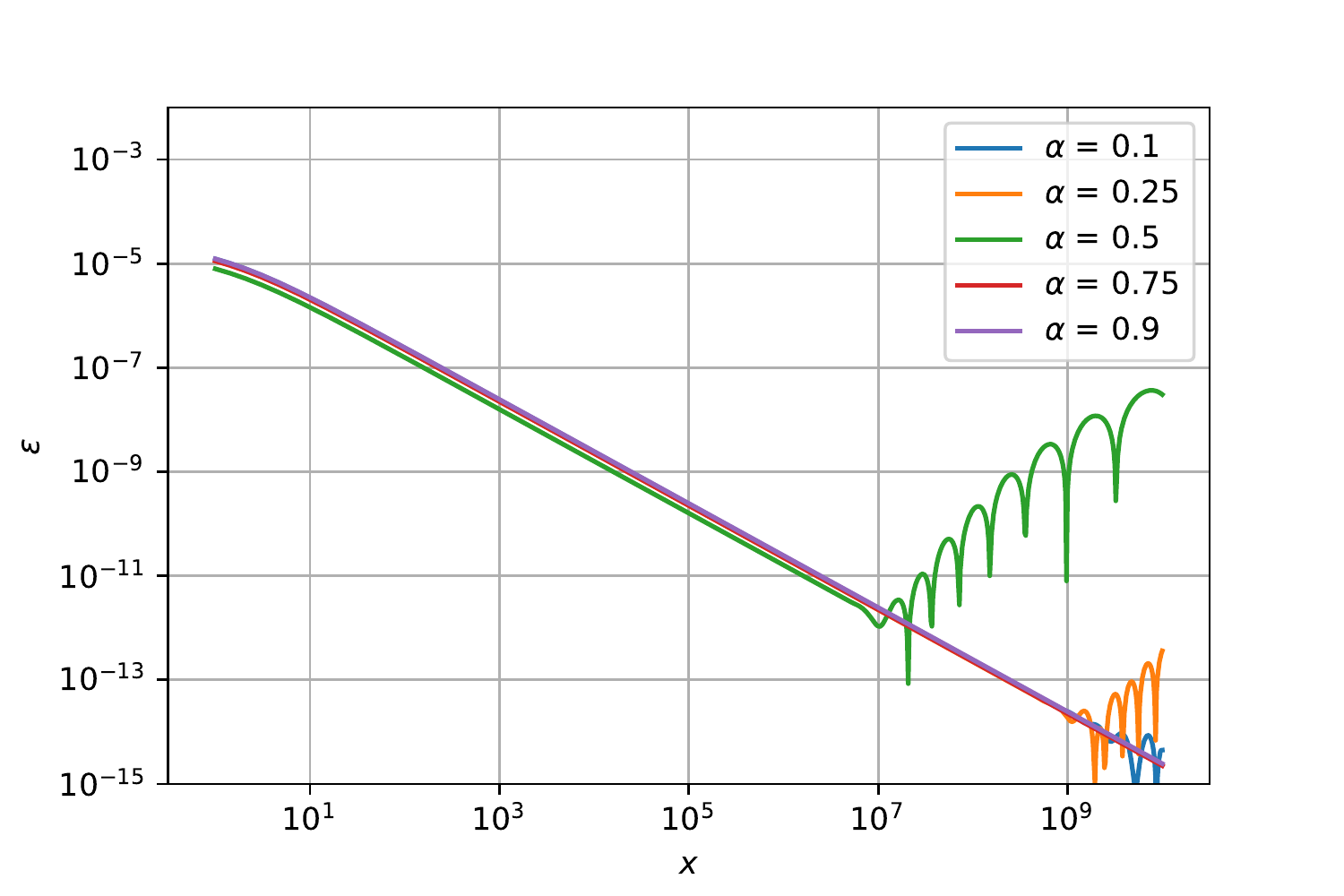}\\
$\varkappa = 3$ \\
\includegraphics[width=\linewidth]{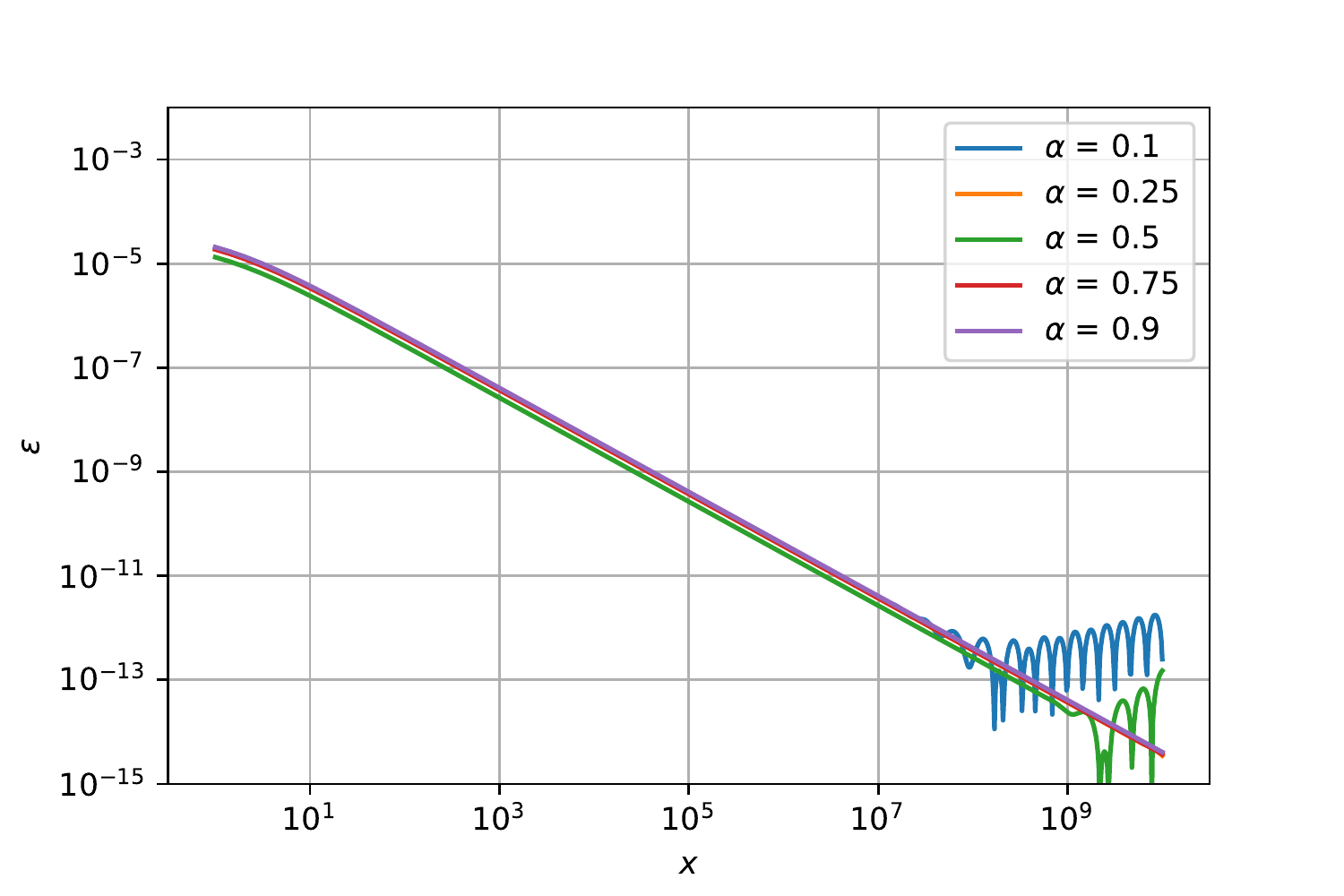}\\
$\varkappa = 5$ \\
\end{minipage}
\begin{minipage}{0.43\linewidth}
\centering
\includegraphics[width=\linewidth]{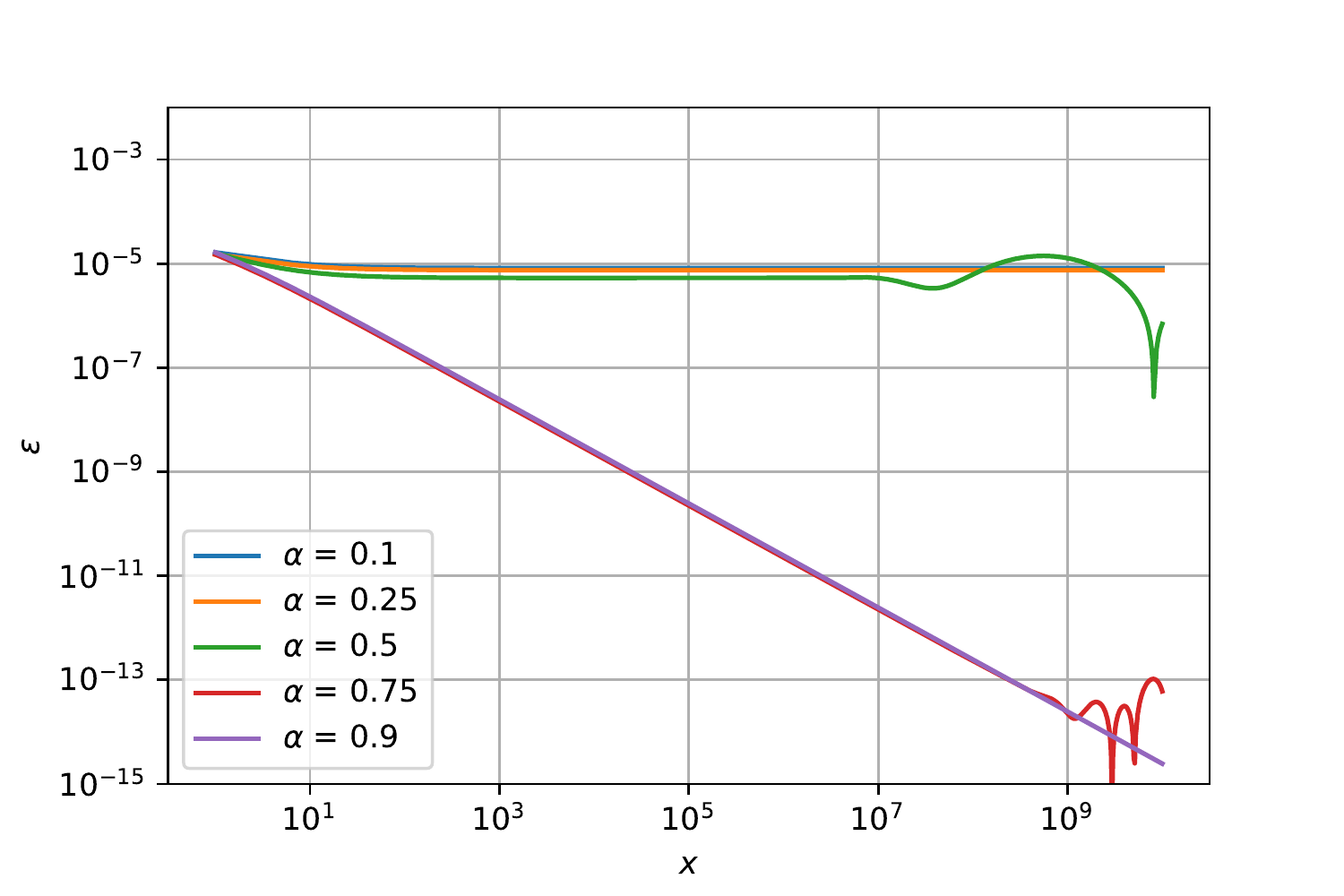}\\
$\varkappa = 2$ \\
\includegraphics[width=\linewidth]{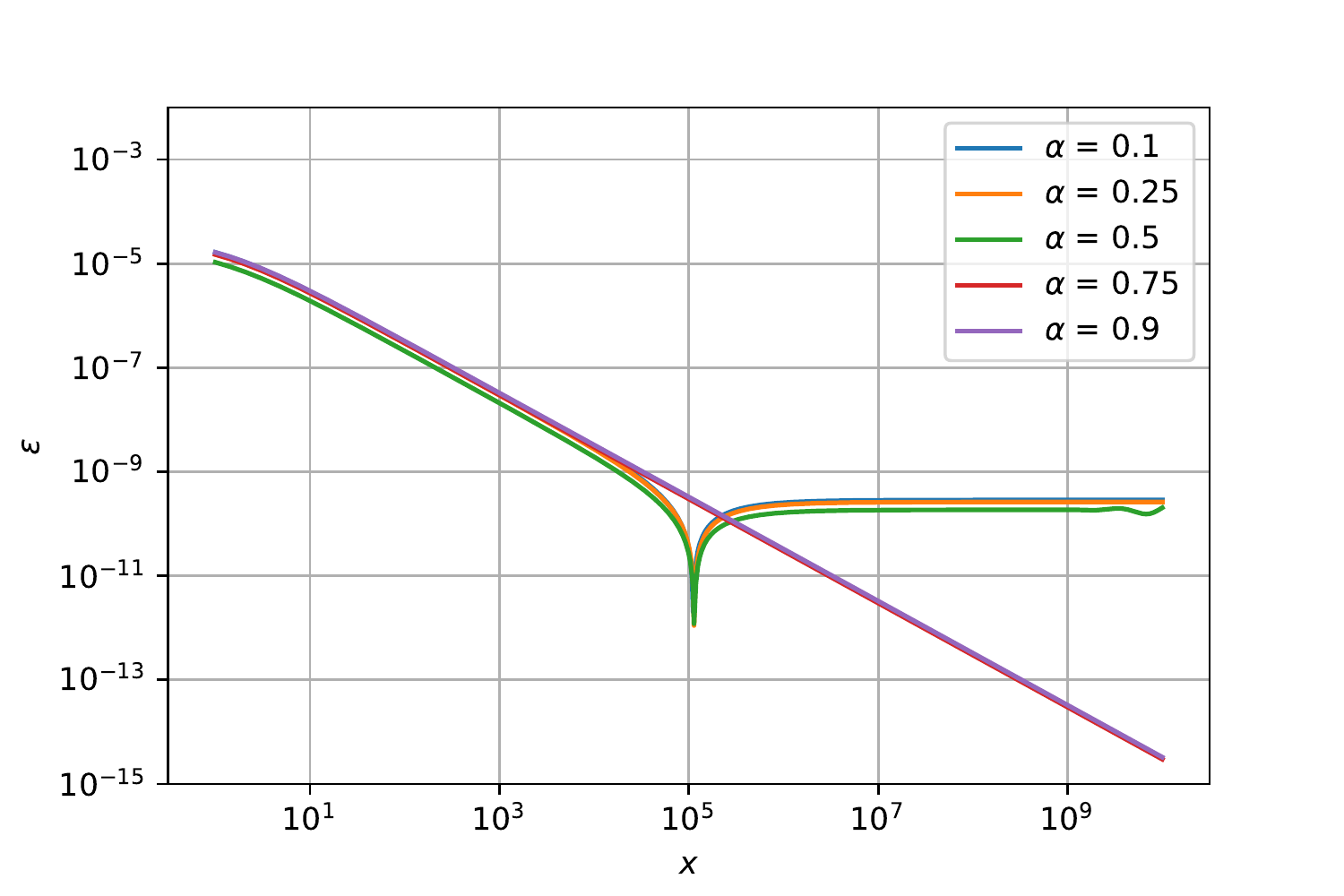}\\
$\varkappa = 4$ \\
\includegraphics[width=\linewidth]{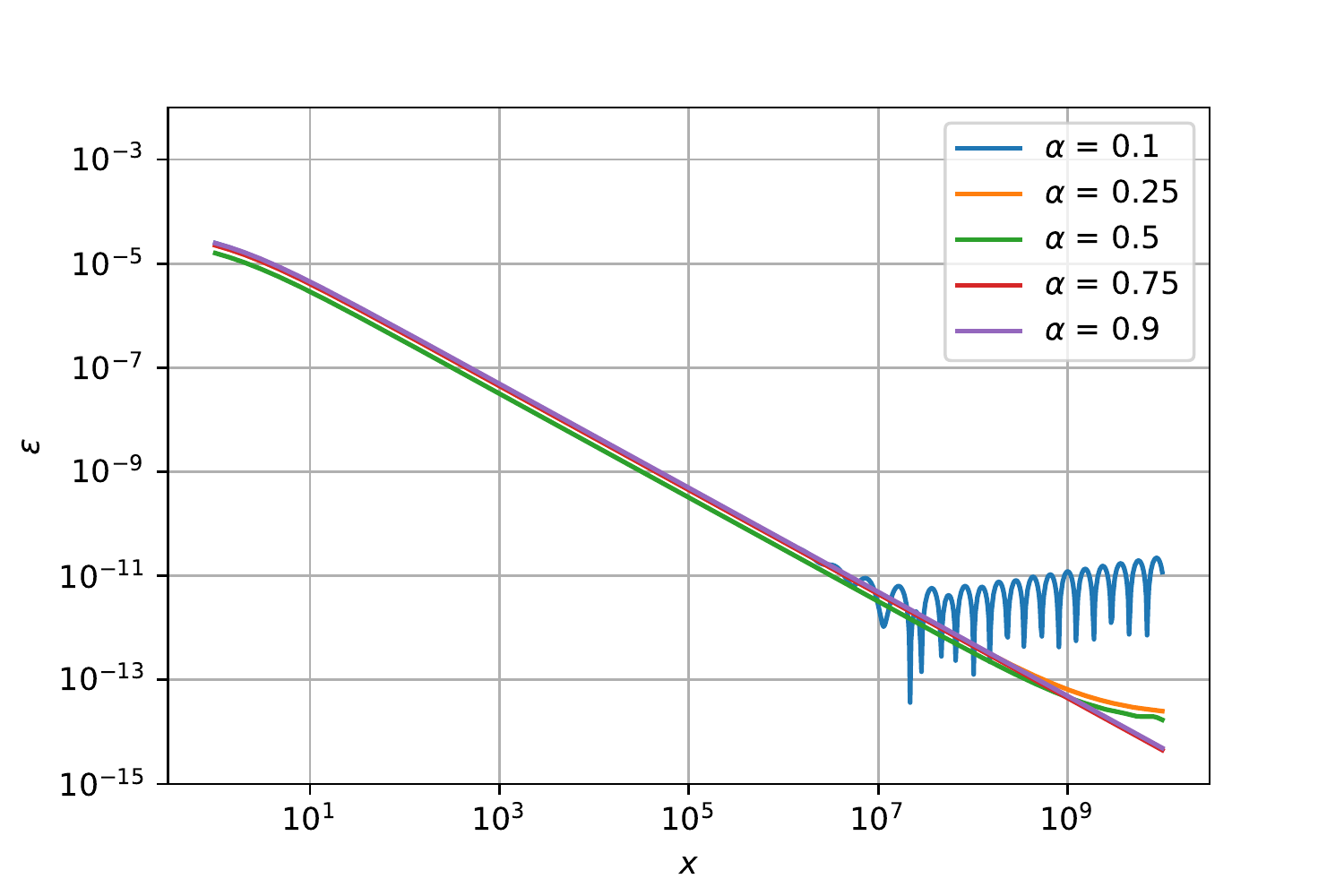}\\
$\varkappa = 6$ \\
\end{minipage}
\caption{Approximation error: the rectangle rule for (\ref{23}), $M = 100$.}
\label{f-2}
\end{figure}

The approximation accuracy of the function $ x^{-\alpha}$ for the representation (\ref{23})
is shown in Figure~\ref{f-2}. The corresponding numerical data are given in Table~2.
We observe similar accuracy results for the representations (\ref{22}) and (\ref{23}).
When focusing on approximations of the fractional power operator, the choice is in favor of using the representation (\ref{22}),
since both of the representations demonstrate comparable accuracy, whereas the computational complexity of the representation (\ref{23}) 
is approximately two times higher.

\begin{center}
\begin{table}[htp]
\label{t-2}
\caption{Approximation error $\overline{\varepsilon}_e$: the rectangle rule for (\ref{23})}
\centering
\begin{tabular}{ccccccc}
\hline
   $M$   &    $\varkappa$   &   $\alpha = 0.1$         &     $\alpha = 0.25$   &    $\alpha = 0.5$   &     $\alpha = 0.75$      &    $\alpha = 0.9$   \\
\hline
        &      1   &  5.347988e-03     &  3.966817e-03     &    1.772738e-03    &  1.500737e-05    &  1.640307e-05  \\
        &      2   &  8.982157e-05     &  1.536406e-04     &    6.366643e-05    &  6.004424e-05    &  6.565671e-05  \\
   50   &      3   &  4.946450e-05     &  4.509810e-05     &    3.184883e-05    &  4.509810e-05    &  4.946450e-05  \\
        &      4   &  6.626068e-05     &  6.021592e-05     &    4.248002e-05    &  6.021592e-05    &  6.626068e-05  \\
        &      5   &  8.336423e-05     &  7.541989e-05     &    5.314105e-05    &  7.541989e-05    &  8.336423e-05  \\
        &      6   &  1.008694e-04     &  9.071839e-05     &    6.381866e-05    &  9.071839e-05    &  1.008694e-04  \\

\hline
        &      1   &  1.289544e-03     &  1.983358e-03     &    8.863758e-04    &  3.751449e-06    &  4.099039e-06  \\
        &      2   &  1.639890e-05     &  3.840905e-05     &    1.591577e-05    &  1.500672e-05    &  1.639890e-05  \\
  100   &      3   &  1.231272e-05     &  1.125908e-05     &    7.958861e-06    &  1.125908e-05    &  1.231272e-05  \\
        &      4   &  1.643508e-05     &  1.501738e-05     &    1.061274e-05    &  1.501738e-05    &  1.643508e-05  \\
        &      5   &  2.057415e-05     &  1.878095e-05     &    1.326849e-05    &  1.878095e-05    &  2.057415e-05  \\
        &      6   &  2.473359e-05     &  2.255024e-05     &    1.592525e-05    &  2.255024e-05    &  2.473359e-05  \\
\hline
        &      1   &  1.638621e-04     &  9.916965e-04     &    4.431897e-04    &  9.378377e-07    &  1.024652e-06  \\
        &      2   &  4.098779e-06     &  9.602172e-06     &    3.978891e-06    &  3.751408e-06    &  4.098779e-06  \\
  200   &      3   &  3.074926e-06     &  2.813809e-06     &    1.989506e-06    &  2.813809e-06    &  3.074926e-06  \\
        &      4   &  4.101019e-06     &  3.752073e-06     &    2.652733e-06    &  3.752073e-06    &  4.101019e-06  \\
        &      5   &  5.128125e-06     &  4.690667e-06     &    3.316076e-06    &  4.690667e-06    &  5.128125e-06  \\
        &      6   &  6.156442e-06     &  5.629614e-06     &    3.979483e-06    &  5.629614e-06    &  6.156442e-06  \\
\hline
\end{tabular}
\end{table}
\end{center}

\begin{figure}
\centering
\begin{minipage}{0.43\linewidth}
\centering
\includegraphics[width=\linewidth]{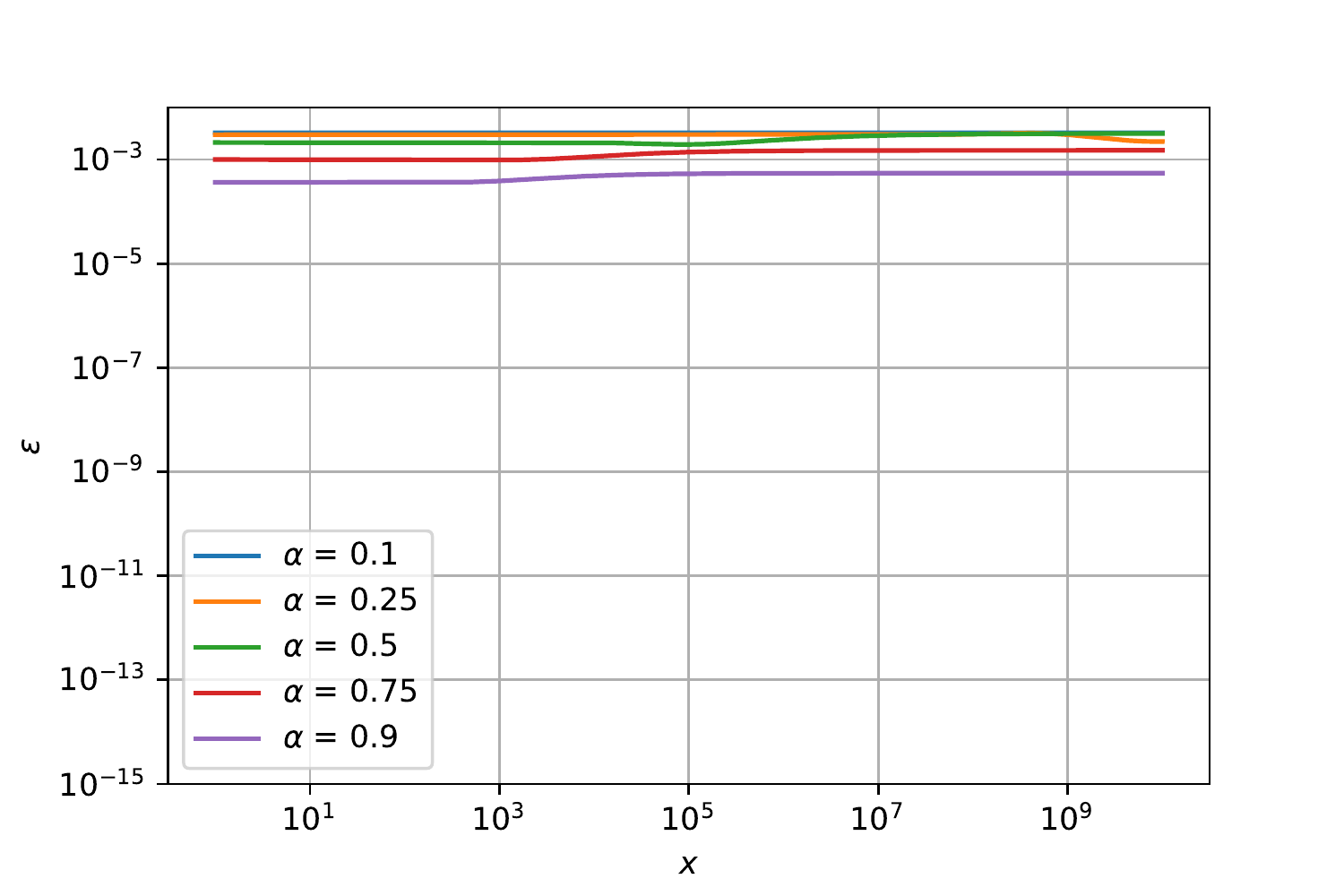}\\
$\varkappa = 1$ \\
\includegraphics[width=\linewidth]{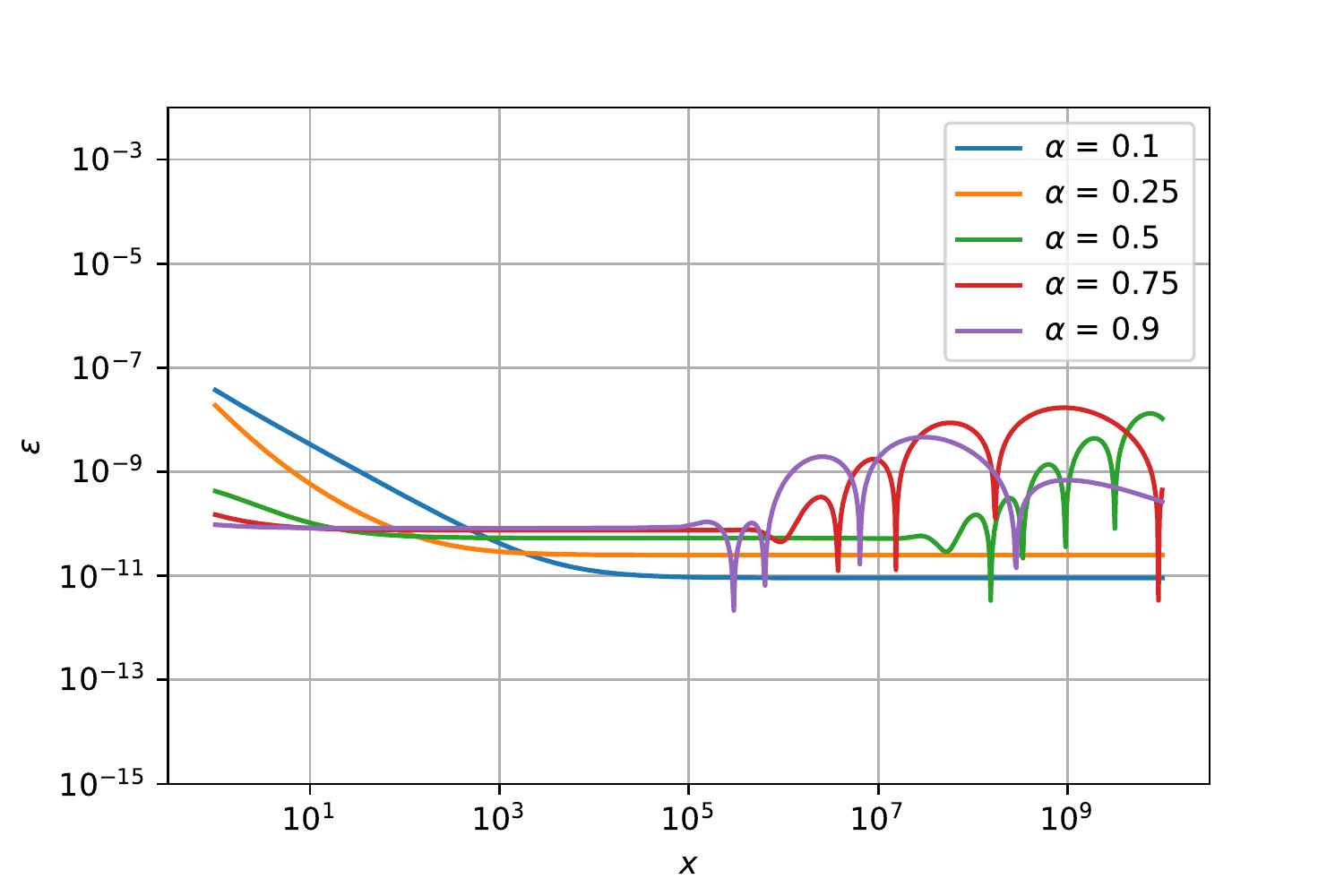}\\
$\varkappa = 3$ \\
\includegraphics[width=\linewidth]{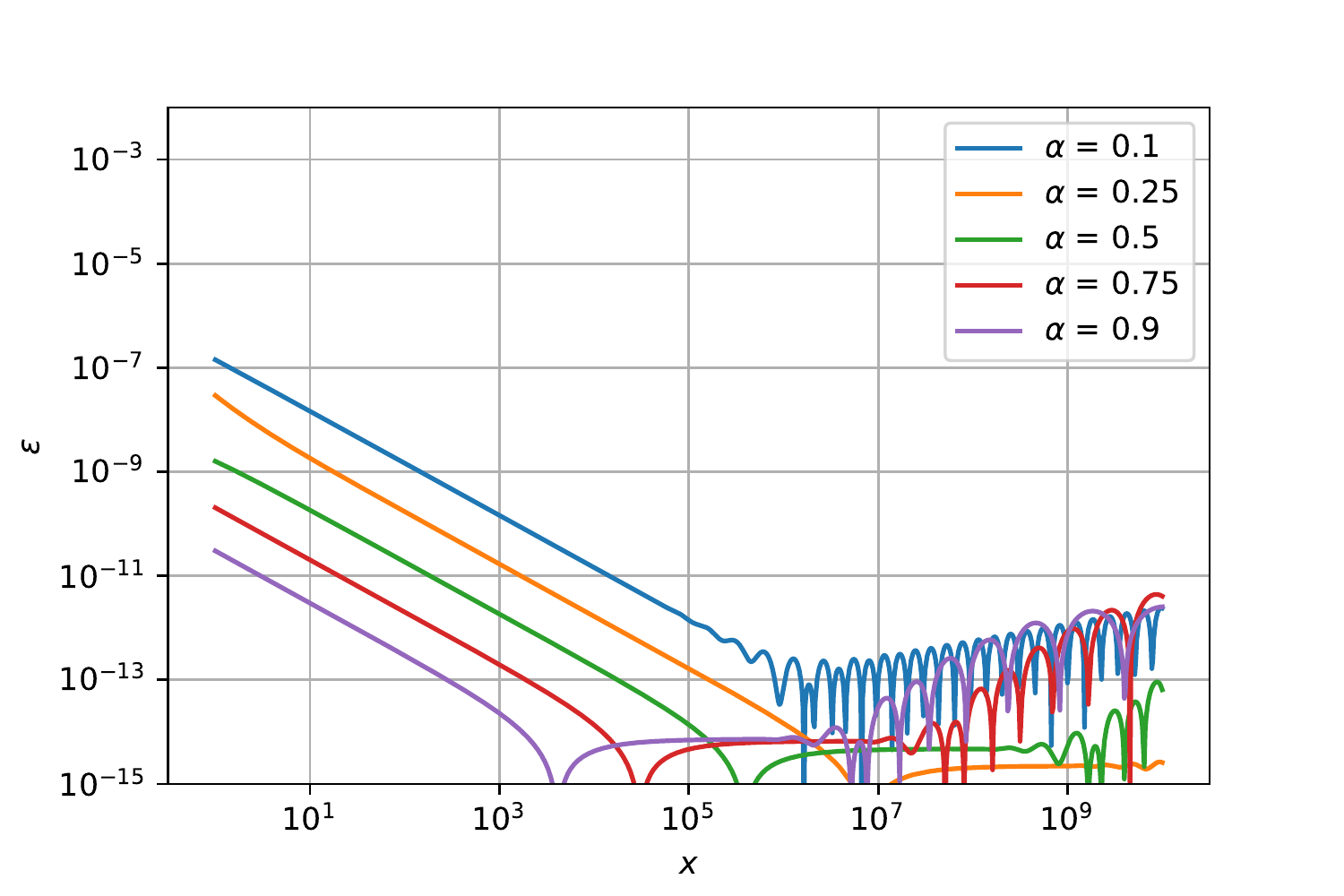}\\
$\varkappa = 5$ \\
\end{minipage}
\begin{minipage}{0.43\linewidth}
\centering
\includegraphics[width=\linewidth]{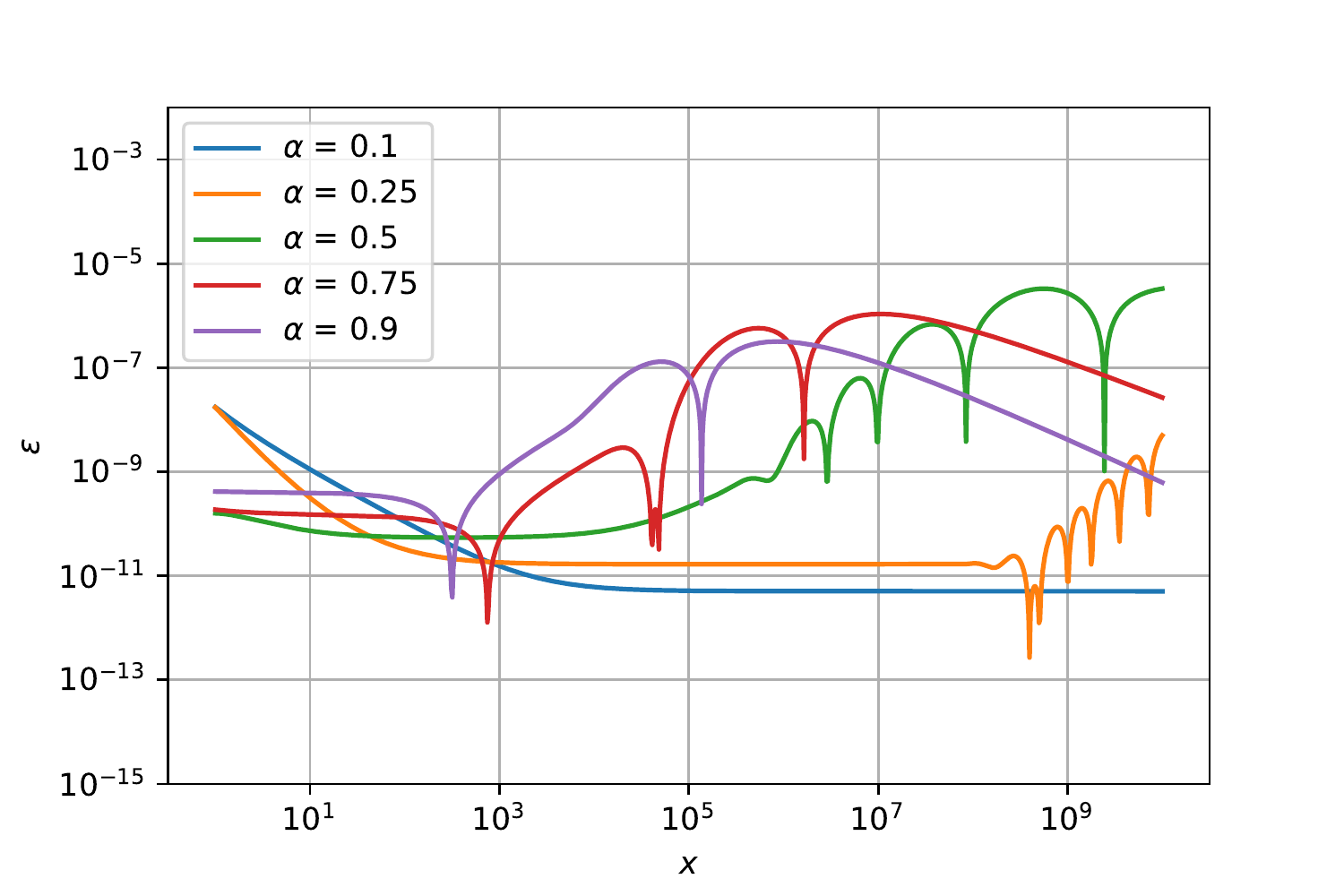}\\
$\varkappa = 2$ \\
\includegraphics[width=\linewidth]{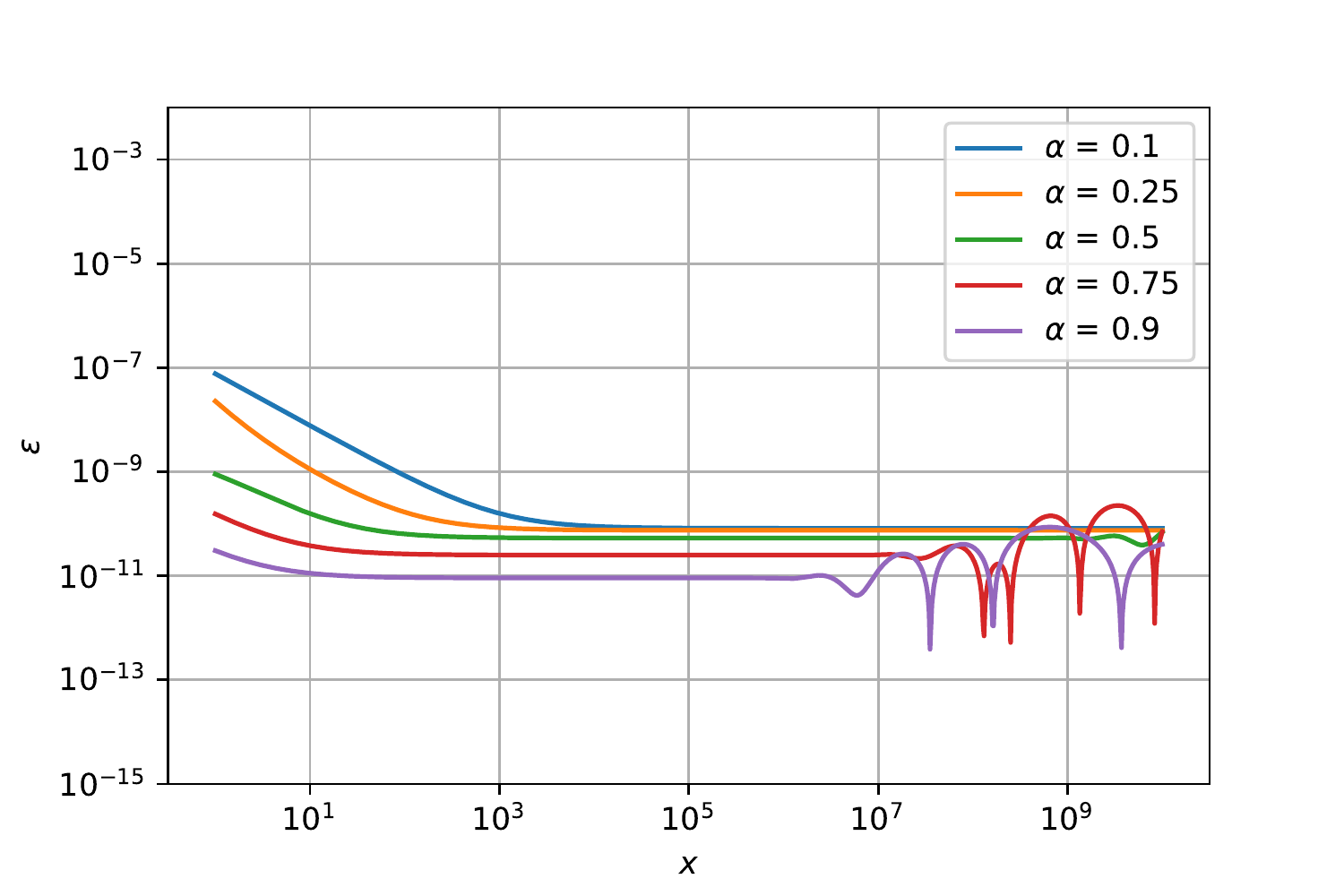}\\
$\varkappa = 4$ \\
\includegraphics[width=\linewidth]{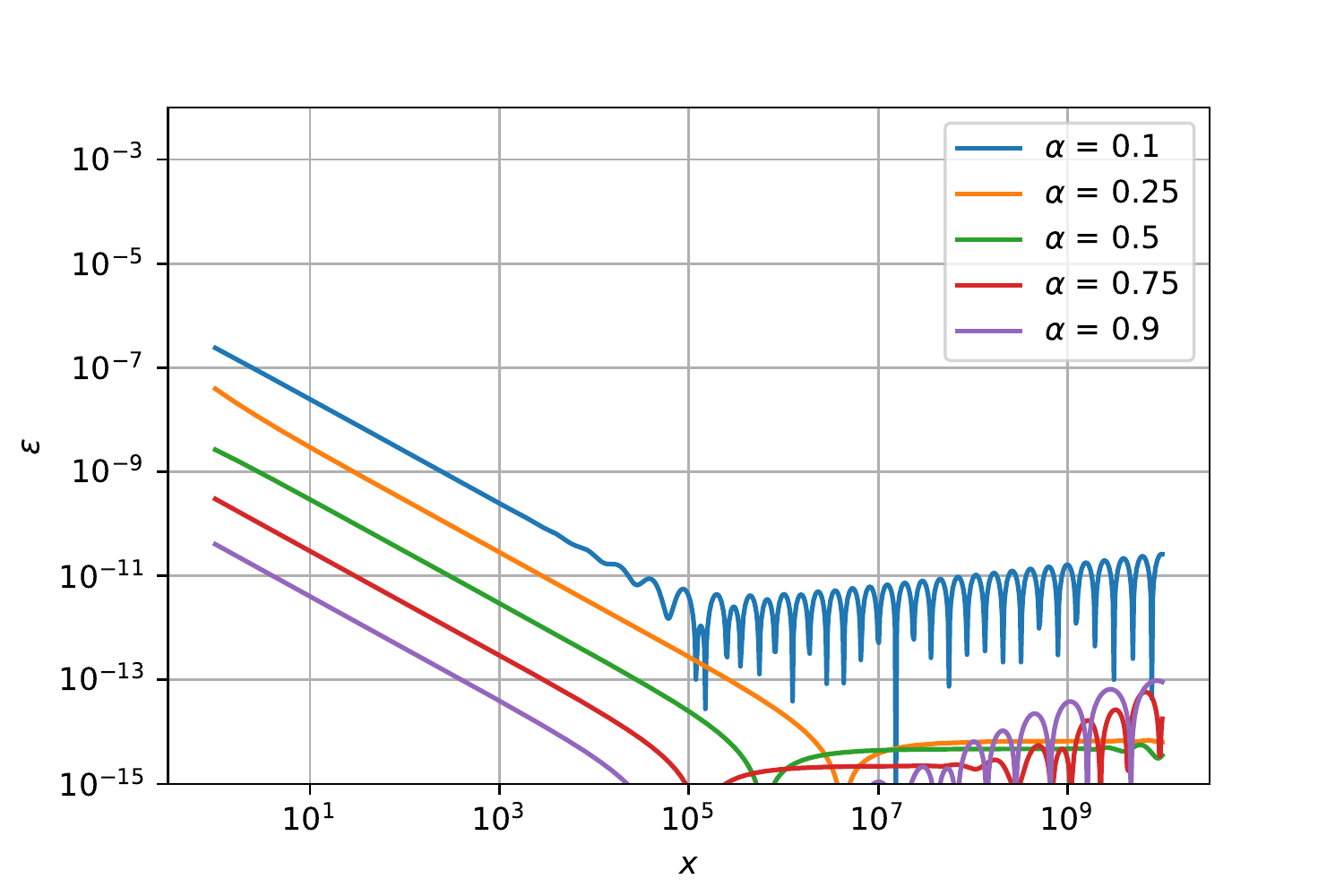}\\
$\varkappa = 6$ \\
\end{minipage}
\caption{Approximation error: Simpson's rule for (\ref{22}), $M = 200$.}
\label{f-3}
\end{figure}

We also present results on the approximation accuracy for Simpson's quadrature rule.
The approximation error of the function $ x^{-\alpha}$ for the representation (\ref{22})
is given in Table~3. 
The theoretical dependence of the accuracy on the number of partitions is observed for large $M$ and 
the parameter $\varkappa \approx 5$. Figure~\ref{f-3} shows the error for a fairly detailed partitioning with $M = 200$.
It is easy to see that the accuracy increases with increasing $\alpha$.

\begin{center}
\begin{table}[htp]
\label{t-3}
\caption{Approximation error $\overline{\varepsilon}_e$: Simpson's rule for (\ref{22})}
\centering
\begin{tabular}{ccccccc}
\hline
   $M$   &    $\varkappa$   &   $\alpha = 0.1$         &     $\alpha = 0.25$   &    $\alpha = 0.5$   &     $\alpha = 0.75$      &    $\alpha = 0.9$   \\
\hline
        &      1   &  1.853227e-02     &  1.799633e-02     &    1.273240e-02    &  6.002109e-03    &  2.185848e-03  \\
        &      2   &  3.949860e-04     &  2.269479e-04     &    5.687915e-05    &  1.750531e-05    &  5.280150e-06  \\
   50   &      3   &  1.275793e-04     &  2.277837e-05     &    4.469762e-06    &  1.174967e-06    &  3.586255e-07  \\
        &      4   &  9.400469e-05     &  5.571588e-06     &    5.614120e-07    &  9.194800e-08    &  3.276035e-08  \\
        &      5   &  2.658768e-04     &  4.032756e-06     &    4.091810e-07    &  5.210468e-08    &  3.395709e-08  \\
        &      6   &  9.293393e-04     &  7.461337e-06     &    6.803497e-07    &  7.645749e-08    &  1.132277e-07  \\
\hline
        &      1   &  9.253872e-03     &  8.993163e-03     &    6.366198e-03    &  3.001054e-03    &  1.092924e-03  \\
        &      2   &  2.695390e-05     &  5.624706e-05     &    1.413623e-05    &  4.328911e-06    &  1.290228e-06  \\
  100   &      3   &  3.871962e-06     &  2.782215e-06     &    5.489516e-07    &  1.401037e-07    &  3.949594e-08  \\
        &      4   &  1.981503e-06     &  3.273869e-07     &    3.283343e-08    &  5.249440e-09    &  1.580743e-09  \\
        &      5   &  2.380378e-06     &  2.949317e-07     &    2.549296e-08    &  3.259146e-09    &  4.821441e-10  \\
        &      6   &  4.575033e-06     &  4.604162e-07     &    4.247676e-08    &  4.799821e-09    &  6.473959e-10  \\
\hline
        &      1   &  3.535332e-03     &  4.491582e-03     &    3.183099e-03    &  1.500527e-03    &  5.464620e-04  \\
        &      2   &  1.864876e-07     &  1.399930e-05     &    3.523765e-06    &  1.076370e-06    &  3.189672e-07  \\
  200   &      3   &  3.715126e-08     &  2.692643e-07     &    6.802848e-08    &  1.710196e-08    &  4.636237e-09  \\
        &      4   &  7.661090e-08     &  2.266696e-08     &    1.986647e-09    &  3.136784e-10    &  8.664345e-11  \\
        &      5   &  1.420196e-07     &  2.919238e-08     &    1.591996e-09    &  2.037560e-10    &  3.016176e-11  \\
        &      6   &  2.398536e-07     &  3.930730e-08     &    2.653151e-09    &  3.000763e-10    &  4.047396e-11  \\
\hline
\end{tabular}
\end{table}
\end{center}

\section{Numerical solution of problems with the fractional power elliptic operator}\label{sec:5} 

Now we apply the quadrature formulas of the rectangle and Simpson's rule considered above for (\ref{22})
to the approximation of the integral representation of the fractional power grid elliptic operator (\ref{16})
for solving the problem (\ref{7}).

\begin{figure}
\centering
\begin{minipage}{0.43\linewidth}
\centering
\includegraphics[width=\linewidth]{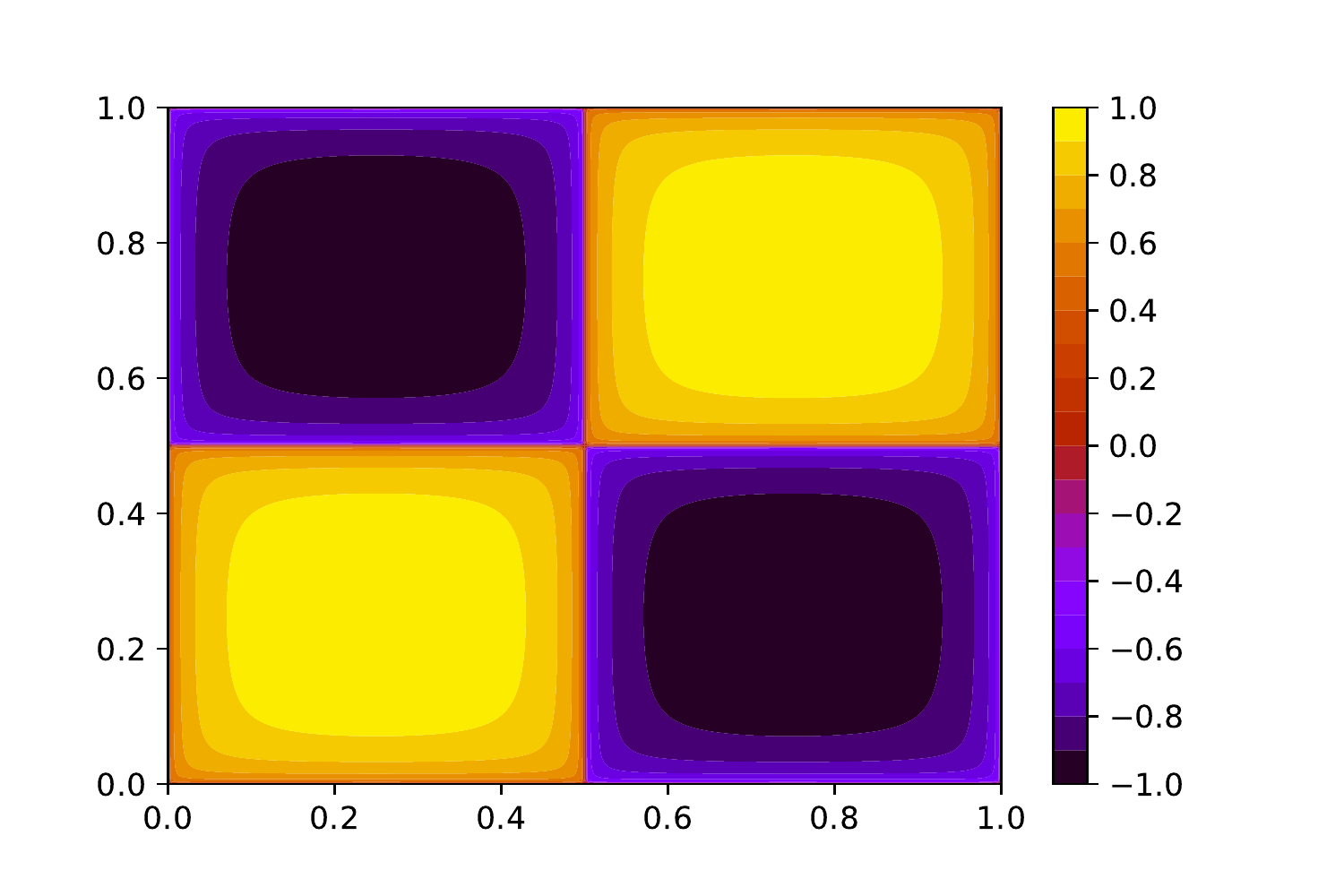}\\
$\alpha  = 0.1, \quad \max u(\bm x) = 6.914610e-01$ \\
\includegraphics[width=\linewidth]{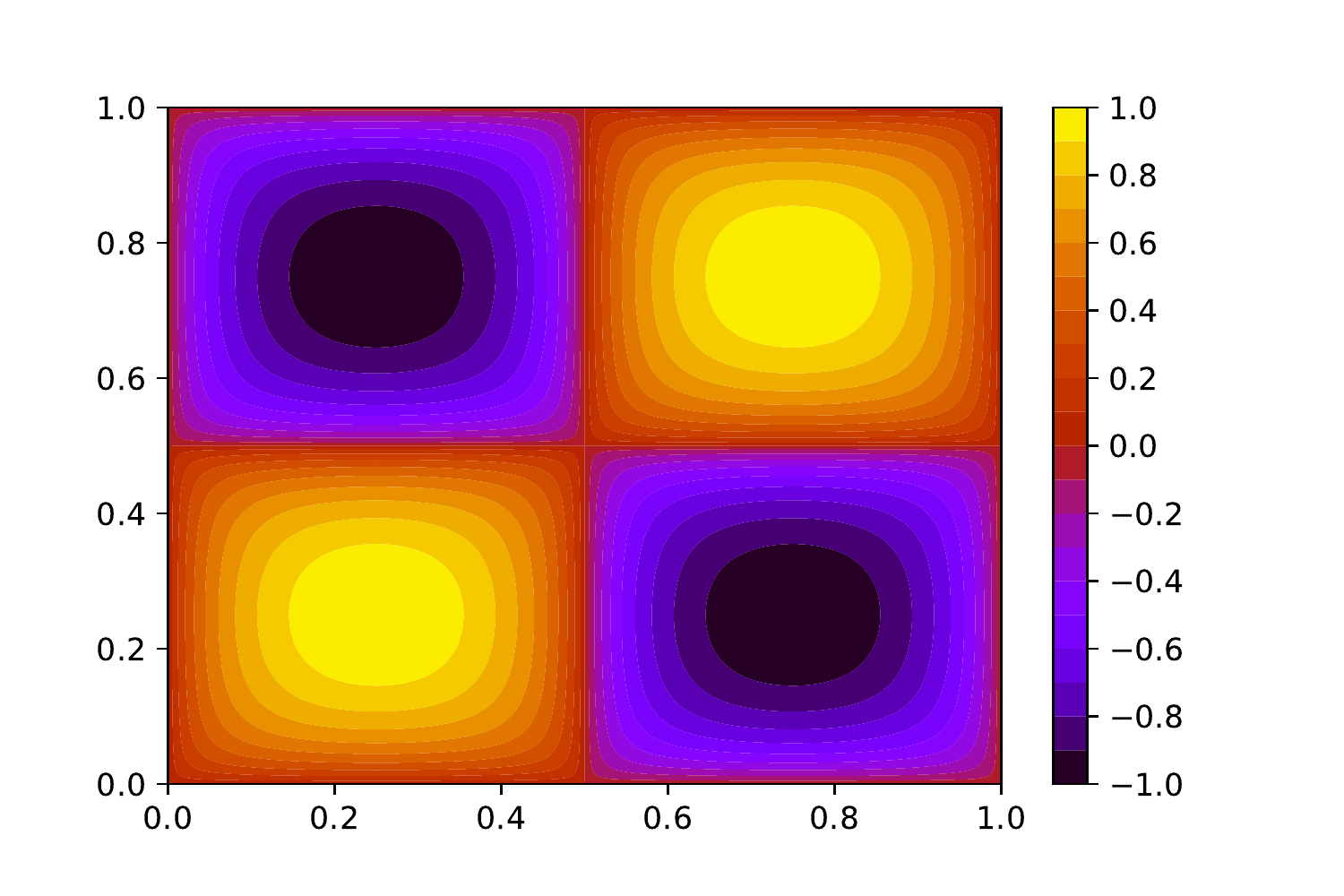}\\
$\alpha  = 0.5, \quad \max u(\bm x) = 1.451668e-01$ \\
\end{minipage}
\begin{minipage}{0.43\linewidth}
\centering
\includegraphics[width=\linewidth]{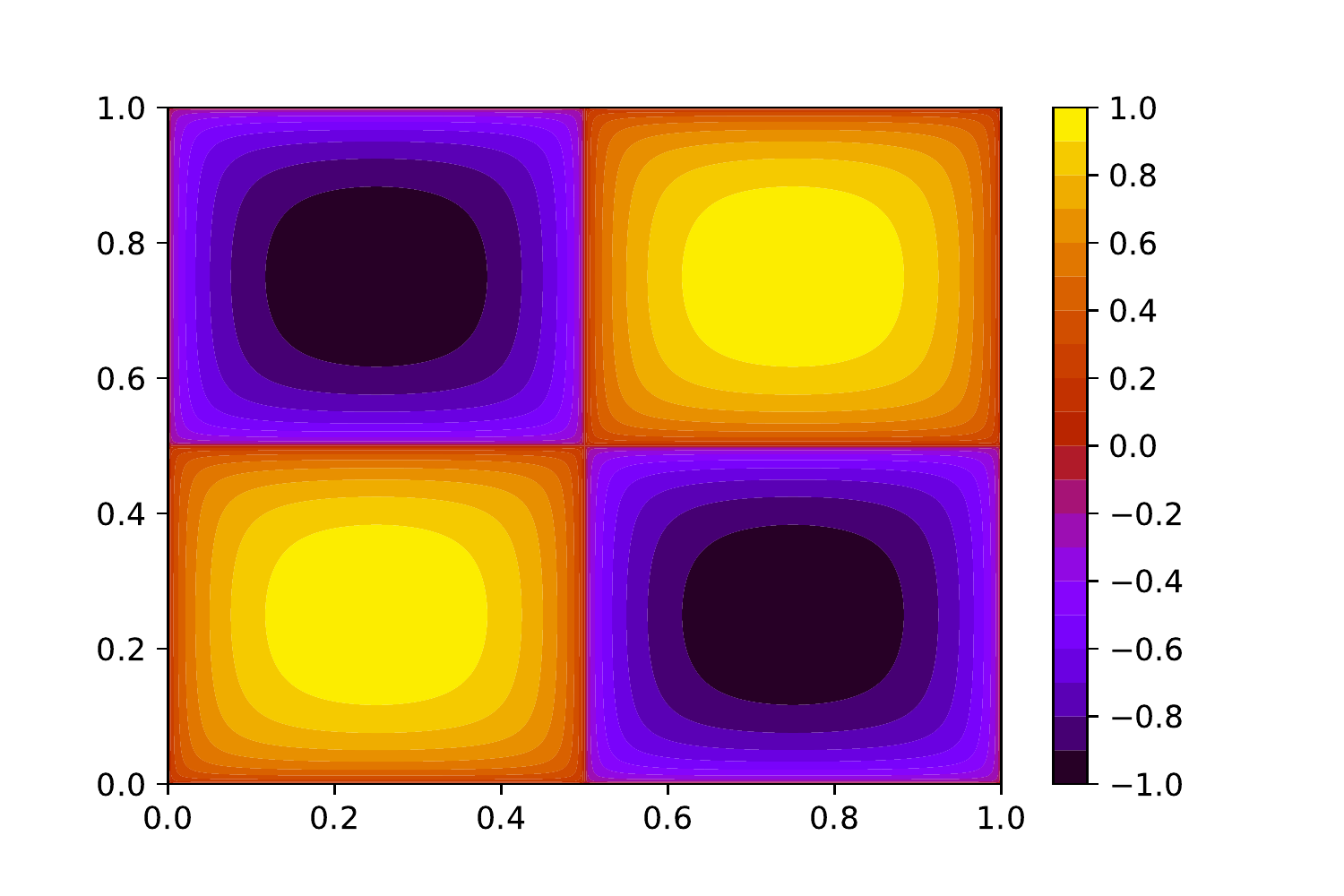}\\
$\alpha  = 0.25, \quad \max u(\bm x) = 3.903881e-01$ \\
\includegraphics[width=\linewidth]{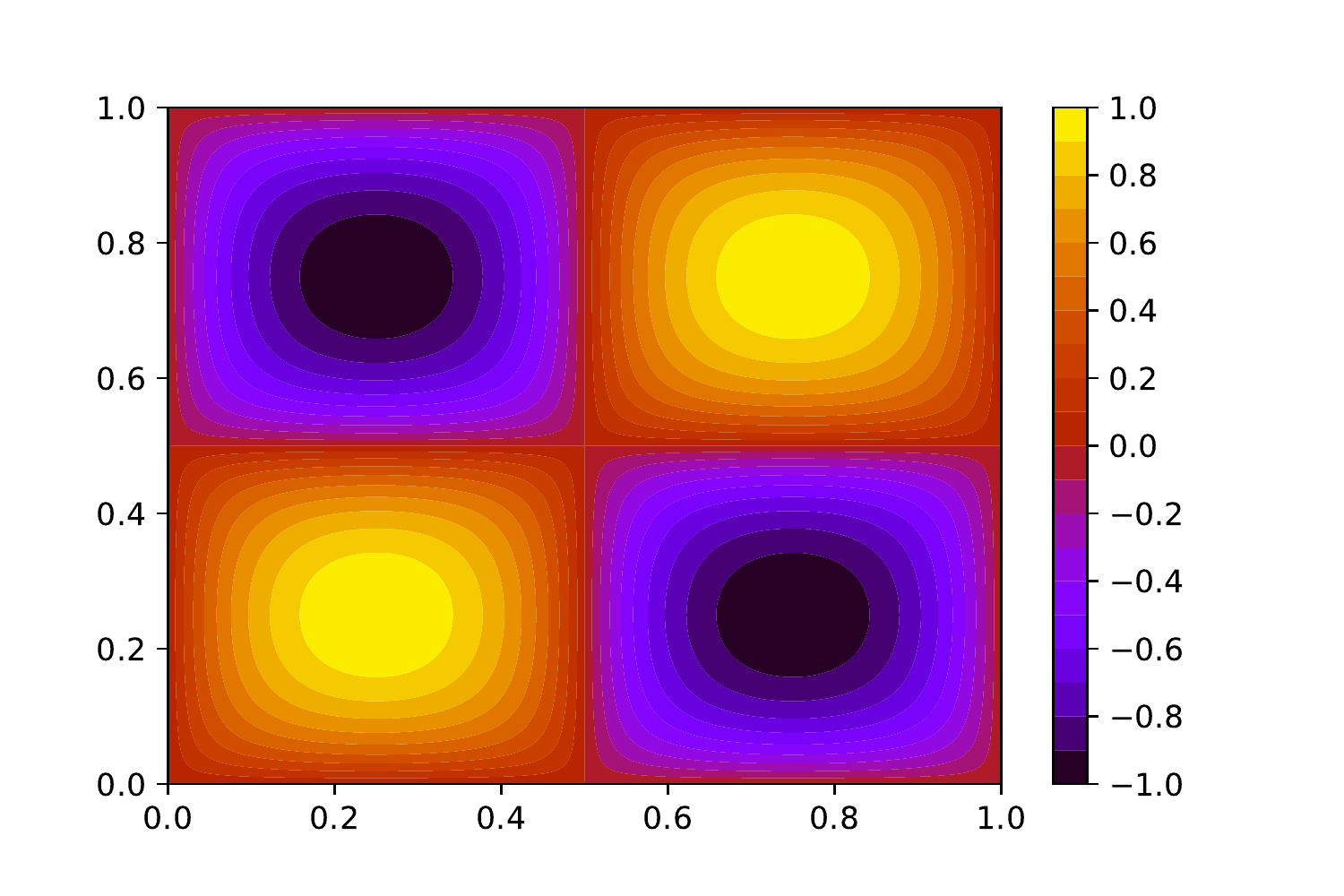}\\
$\alpha  = 0.75, \quad \max u(\bm x) = 5.227385e-02$ \\
\end{minipage}
\caption{Normalized solution of the problem (\ref{3}), (\ref{4}), (\ref{6}), (\ref{24}) on the grid $N_1 = N_2 = 256$.}
\label{f-4}
\end{figure}

In the results presented below, the calculations were performed on the grid with $N_1 = N_2 = 256$.
The peculiarities of boundary value problems with the fractional power operator were studied on the problem
(\ref{3}), (\ref{4}), (\ref{6}) in the unit square ($l_1 = l_2 = 1$) with the right-hand side
\begin{equation}\label{24}
 f(\bm x) = \mathrm{sgn}(x_1 - 0.5) \ \mathrm{sgn}(x_2-0.5),
 \quad  \mathrm{sgn}(x) := \left \{
 \begin{array}{ccc}
  -1,  &  x < 0, \\
   0,  &  x = 0,  \\
   1,  &  x > 0 .  \\
\end{array}
\right .  
\end{equation}
The numerical solution of the problem (\ref{7}) with the right-hand side (\ref{24})
is depicted in Figure~\ref{f-4} for different values of $\alpha$.
For convenience of the comparison, there is shown the function
\[
 y(\bm x) = \frac{1}{\max u(\bm x)} u(\bm x) ,
 \quad \bm x \in \omega . 
\]
For the discontinuous right-hand side (\ref{24}),
we observe the formation of internal boundary layers with decreasing $\alpha$.
 
The accuracy of the approximate solution of the problem with the fractional power of the Laplace operator obtained using
the rectangle rule is illustrated by the data in Table~4.
Here the relative error is determined as follows:
\[
 \varepsilon = \frac{\|w-u\|_2}{\|u\|_2} ,
 \quad \|u\|_2^2 = \sum_{\bm x \in \omega} u^2(\bm x) h_1 h_2 ,
 \quad \varepsilon_\infty  = \frac{\|w-u\|_\infty }{\|u\|_\infty } ,
 \quad \|u\|_\infty  = \max_{\bm x \in \omega} |u(\bm x)| ,
\] 
where $u$ is the exact solution of the problem (\ref{7}), (\ref{24}) and $w$ is the approximate one.
Similar data for Simpson's quadrature formula are shown in Table~5.
In these calculations, we employed $\varkappa = 3$ in the representation (\ref{22}) for the rectangle rule
and $\varkappa = 5$ for Simpson's rule.
We observe good accuracy of the approximate solution even for a small number of nodes of the quadrature formula.

\begin{center}
\begin{table}[htp]
\label{t-4}
\caption{Solution error for the problem (\ref{7}), (\ref{24}): the rectangle rule}
\centering
\begin{tabular}{ccccccc}
\hline
   $M$   &    error  &   $\alpha = 0.1$         &     $\alpha = 0.25$   &    $\alpha = 0.5$   &     $\alpha = 0.75$      &    $\alpha = 0.9$   \\
\hline
   50   &  $\varepsilon$  &  1.669630e-06     &  3.172818e-06     &    6.976081e-06    &  9.886524e-06    &  6.807410e-06  \\
        &  $\varepsilon_{\infty}$  &  1.931587e-06     &  3.483513e-06     &    7.301196e-06    &  1.006909e-05    &  6.868003e-06  \\
\hline
  100   &  $\varepsilon$  &  4.234571e-07     &  7.955079e-07     &    1.749686e-06    &  2.493483e-06    &  1.746540e-06  \\
        &  $\varepsilon_{\infty}$  &  4.899010e-07     &  8.733658e-07     &    1.830969e-06    &  2.538547e-06    &  1.758218e-06  \\
\hline
  200   &  $\varepsilon$  &  1.062490e-07     &  1.990234e-07     &    4.377758e-07    &  6.247408e-07    &  4.394942e-07  \\
        &  $\varepsilon_{\infty}$  &  1.229199e-07     &  2.184997e-07     &    4.580968e-07    &  6.359638e-07    &  4.422973e-07  \\
\hline
\end{tabular}
\end{table}
\end{center}

\begin{center}
\begin{table}[htp]
\label{t-5}
\caption{Solution error for the problem (\ref{7}), (\ref{24}): Simpson's rule}
\centering
\begin{tabular}{ccccccc}
\hline
   $M$   &    error  &   $\alpha = 0.1$         &     $\alpha = 0.25$   &    $\alpha = 0.5$   &     $\alpha = 0.75$      &    $\alpha = 0.9$   \\
\hline
   50   &  $\varepsilon$  &  1.286770e-04     &  2.641118e-07     &    5.182607e-08    &  2.110766e-08    &  2.297871e-07  \\
        &  $\varepsilon_{\infty}$  &  1.558345e-04     &  3.437617e-07     &    5.442683e-08    &  3.023970e-08    &  2.475985e-07  \\
\hline
  100   &  $\varepsilon$  &  7.460526e-08     &  9.578675e-09     &    3.256882e-09    &  1.074501e-09    &  2.876888e-10  \\
        &  $\varepsilon_{\infty}$  &  1.013479e-07     &  1.051652e-08     &    3.407962e-09    &  1.094340e-09    &  2.910291e-10  \\
\hline
  200   &  $\varepsilon$  &  2.528173e-09     &  6.061693e-10     &    2.038772e-10    &  6.772998e-11    &  1.906629e-11  \\
        &  $\varepsilon_{\infty}$  &  2.924801e-09     &  6.656861e-10     &    2.133260e-10    &  6.895408e-11    &  1.920284e-11  \\
\hline
\end{tabular}
\end{table}
\end{center}

For problems with the fractional power elliptic operator, we considered separately
the dependence of the accuracy of the approximate solution on the smoothness of  the exact solution.
Above, we presented the results for the case of the discontinuous right-hand side (see (\ref{24})).
For
\begin{equation}\label{25}
 f(\bm x) = x_1 x_2 ,
\end{equation} 
we have a boundary layer on a part of the boundary of the computational domain for small $\alpha$.
The results on the accuracy of the approximate solution for this case are given in Table~6.
Similar data for the smoother exact solution obtained for
\begin{equation}\label{26}
 f(\bm x) = x_1(1-x_1) x_2 (1-x_2) ,
\end{equation}
are presented in Table~7.
The numerical results indicate that the accuracy of the solution of the discrete problem (\ref{7}) is practically independent
from the smoothness of the solution of the continuous problem (\ref{6}).

\begin{center}
\begin{table}[htp]
\label{t-6}
\caption{Solution error for the problem (\ref{7}), (\ref{25}): Simpson's rule}
\centering
\begin{tabular}{ccccccc}
\hline
   $M$   &    error  &   $\alpha = 0.1$         &     $\alpha = 0.25$   &    $\alpha = 0.5$   &     $\alpha = 0.75$      &    $\alpha = 0.9$   \\
\hline
   50   &  $\varepsilon$  &  1.209581e-04     &  3.444934e-07     &    9.686650e-08    &  2.433874e-08    &  7.815106e-08  \\
        &  $\varepsilon_{\infty}$  &  1.141529e-04     &  3.688685e-07     &    9.631475e-08    &  2.422181e-08    &  8.320527e-08  \\
\hline
  100   &  $\varepsilon$  &  1.030346e-07     &  2.380002e-08     &    6.073046e-09    &  1.493598e-09    &  3.465175e-10  \\
        &  $\varepsilon_{\infty}$  &  1.108808e-07     &  2.349288e-08     &    6.059424e-09    &  1.493626e-09    &  3.468453e-10  \\
\hline
  200   &  $\varepsilon$  &  7.245613e-09     &  1.539020e-09     &    3.799707e-10    &  9.359799e-11    &  2.202853e-11  \\
        &  $\varepsilon_{\infty}$  &  6.960763e-09     &  1.518666e-09     &    3.791153e-10    &  9.359548e-11    &  2.203251e-11  \\
\hline
\end{tabular}
\end{table}
\end{center}

\begin{center}
\begin{table}[htp]
\label{t-7}
\caption{Solution error for the problem (\ref{7}), (\ref{26}): Simpson's rule}
\centering
\begin{tabular}{ccccccc}
\hline
   $M$   &    error  &   $\alpha = 0.1$         &     $\alpha = 0.25$   &    $\alpha = 0.5$   &     $\alpha = 0.75$      &    $\alpha = 0.9$   \\
\hline
   50   &  $\varepsilon$  &  5.800855e-05     &  4.000618e-07     &    1.056566e-07    &  2.455535e-08    &  4.044510e-08  \\
        &  $\varepsilon_{\infty}$  &  6.236828e-05     &  3.995393e-07     &    1.061651e-07    &  2.460929e-08    &  3.914376e-08  \\
\hline
  100   &  $\varepsilon$  &  6.407204e-08     &  2.899259e-08     &    6.627865e-09    &  1.541997e-09    &  3.513618e-10  \\
        &  $\varepsilon_{\infty}$  &  6.549602e-08     &  2.924586e-08     &    6.659433e-09    &  1.545095e-09    &  3.516994e-10  \\
\hline
  200   &  $\varepsilon$  &  9.952338e-09     &  1.877999e-09     &    4.146733e-10    &  9.659714e-11    &  2.227958e-11  \\
        &  $\varepsilon_{\infty}$  &  1.006940e-08     &  1.894607e-09     &    4.166460e-10    &  9.678749e-11    &  2.229965e-11  \\
\hline
\end{tabular}
\end{table}
\end{center}

\section{Conclusions}\label{sec:6} 

\begin{enumerate}
 \item Different variants of the integral representation of the function $x^{-\alpha}, \ x \geq 1, \ 0 < \alpha < 1$
are considered for approximating the fractional power operator. There are highlighted the integral representations,
where the integrand has no singularities.
By introducing new integration variables, integral representations are proposed for the function $x^{-\alpha}$,
where the derivatives of the integrand are bounded.
 \item The rectangle and Simpson quadrature formulas were constructed to approximate the function $x^{-\alpha}$
for various values of $0 < \alpha < 1$.
The accuracy of these quadrature formulas is numerically investigated depending on the key approximation parameters.
 \item The model problem with the fractional power of the Laplace operator is considered in a rectangle.
Its approximate solution is obtained on the basis of the approximation of the function $x^{-\alpha}$ applying the
rectangle and Simpson's quadrature rule. Calculations demonstrate high accuracy of the approximate solution
if we use about one hundred quadrature nodes.
\end{enumerate} 

\section*{Acknowledgements}
The publication 
was financially supported by the Ministry of Education and Science of the Russian Federation (the Agreement number 02.a03.21.0008).


\end{document}